\documentclass[11pt]{article}
\usepackage{amsfonts}
\usepackage{amssymb,latexsym}
\usepackage{graphicx,amsmath,amscd}
\usepackage[pdftitle={Paper}, pdfauthor={Ortiz}]{hyperref}
\usepackage{stmaryrd}
\usepackage{hyperref}  
\usepackage{enumerate}
\usepackage{amssymb}
\usepackage{amsthm}
\usepackage{xspace}
\usepackage{latexsym}
\usepackage{array}
\usepackage[T1]{fontenc}
\usepackage{enumitem}
\usepackage[all]{xy}

\newcommand{\rmap}       {\longrightarrow}
\newcommand{\rr}             {\rightrightarrows}

\newcommand{\Hom}        {{\mathrm {Hom}}}

\newcommand{\G}         {\mathcal{G}}

\newcommand{\la}           {\mathcal{LA}}
\newcommand{\vb}       {\mathcal{VB}}

\newcommand{\Cour}[1]      {[\![#1]\!]}

\theoremstyle{plain}
\newtheorem{theorem}{Theorem}[section]
\newtheorem{proposition}{Proposition}[section]
\newtheorem{lemma}{Lemma}[section]
\newtheorem{corollary}{Corollary}[section]
\theoremstyle{definition}
\newtheorem{example}{Example}[section]
\newtheorem{definition}{Definition}[section]
\newtheorem{remark}{Remark}[section]

\newcommand\R{\mathbb{R}}

\newcommand\N{\mathbb{N}}

\DeclareMathAlphabet{\mathpzc}{OT1}{pcz}{m}{it}

\begin{document}

\author{C. Ortiz\thanks{{\tt cortiz@ime.usp.br}, {Instituto de Matem\'atica e Estat\'istica, Universidade de S\~ao Paulo, Rua do Mat\~ao 1010, Cidade Universit\'aria, 05508-090, S\~ao Paulo - Brasil}}, J. Waldron\thanks{\texttt{james.waldron@ncl.ac.uk}, {School of Mathematics and Statistics, Newcastle University, Newcastle upon Tyne NE1 7RU, UK.}}}

\date{\today}
\title{On the Lie 2-algebra of sections of an $\mathcal{LA}$-groupoid}

\maketitle

\begin{abstract}

In this work we introduce the category of multiplicative sections of an $\la$-groupoid. We prove that these categories carry natural strict Lie 2-algebra structures, which are Morita invariant. As applications, we study the algebraic structure underlying multiplicative vector fields on a Lie groupoid and in particular vector fields on differentiable stacks. We also introduce the notion of geometric vector field on the quotient stack of a Lie groupoid, showing that the space of such vector fields is a Lie algebra. We describe the Lie algebra of geometric vector fields in several cases, including classifying stacks, quotient stacks of regular Lie groupoids and in particular orbifolds, and foliation groupoids.

\end{abstract}

\tableofcontents

\section{Introduction}

The concept of $\la$-groupoid first appeared in the work of K. Mackenzie \cite{mackI} as the infinitesimal counterpart of double Lie groupoids. Intuitively, an $\la$-groupoid is thought of as a groupoid object in the category of Lie algebroids. Besides Lie algebroids of double groupoids, geometric examples of $\la$-groupoids include the cotangent bundle of a Poisson groupoid \cite{mackII} and more generally multiplicative Dirac structures on Lie groupoids \cite{Ortiz13}. 

Additionally, $\la$-groupoids have recently appeared in connection with the notion of Lie algebroid over a differentiable stack, introduced by the second author in \cite{Waldronthesis}. Roughly, differentiable stacks are generalized quotients of smooth manifolds which arise naturally in moduli problems, including orbifolds as particular instances. It is well known that isomorphism classes of differentiable stacks correspond to Morita equivalence classes of Lie groupoids. Hence, one can think of a differentiable stack as a Lie groupoid up to Morita equivalence, e.g. an orbifold corresponds to a proper and \'etale Lie groupoid up to Morita equivalence. The notion of Lie algebroid over a differentiable stack was introduced in \cite{Waldronthesis}, where it is shown that a Lie algebroid over a differentiable stack corresponds to a special type of $\la$-groupoid. Hence, when working with Lie algebroids over a differentiable stack, one can either work with the stacky version or with the $\la$-groupoid picture. In this work we will follow the second approach. The main goal of the paper is to study the algebraic structure underlying the space of sections of an $\la$-groupoid. Our guiding principle is thinking of an $\la$-groupoid as a categorified version of a Lie algebroid, hence the space of sections should be a categorified version of a Lie algebra. This naturally leads us to the notion of Lie 2-algebra in the sense of Baez and Crans \cite{BaezCrans}, i.e.  an internal category to the category of Lie algebras. 

Our approach is as follows. The forgetful functor from the category of Lie algebroids to the category of vector bundles, extends to a forgetful functor from the category of $\la$-groupoids to the category of $\vb$-groupoids. Hence we first study the space of sections of a $\vb$-groupoid and then we restrict our study to $\la$-groupoids. The notion of $\vb$-groupoid was introduced by Pradines in \cite{pradines} in connection with symplectic groupoids. Also, $\vb$-groupoids have been further studied by K. Mackenzie \cite{mackbook} and their relation with representation theory was shown by Mehta and Gracia Saz in \cite{GM10}. Also, $\vb$-groupoids play a role in a possible notion of 2-vector bundle over a differentiable stack \cite{dHO}.

Given a $\vb$-groupoid $\mathcal{V}\rr E$ over $\G\rr M$, we introduce the category of multiplicative sections of $\mathcal{V}$, extending the category of multiplicative vector fields introduced by Hepworth in \cite{hepworth}. The main result concerning multiplicative sections is Theorem \ref{2vect}, which states that the category of multiplicative sections of $\mathcal{V}$ has a natural structure of 2-vector space, induced by a 2-term complex of vector spaces canonically associated to $\mathcal{V}$. Then we move to the case of $\la$-groupoids and we prove Theorem \ref{thm:crossed}, which says that the 2-term complex of vector spaces associated to the underlying $\vb$-groupoid is part of a structure of crossed module of Lie algebras. As a consequence, we show Theorem \ref{Lie2} which states that the category of multiplicative sections of any $\la$-groupoid inherits a canonical structure of Lie 2-algebra. In the special case of the tangent $\la$-groupoid, Theorem \ref{Lie2} gives a positive answer to a conjecture by Hepworth \cite{hepworth} about the Lie 2-algebra structure on the category of multiplicative vector fields on a Lie groupoid. 

In order to obtain a geometric structure in the level of differentiable stacks, we need to deal with Morita invariance of the category of multiplicative sections or equivalently of the corresponding 2-term complexes. Morita maps between $\vb$-groupoids were introduced in \cite{dHO}, and allow one to define a notion of Morita equivalence for $\mathcal{VB}$-groupoids. In Theorem \ref{moritainvarianceVB2} we show that Morita equivalent $\mathcal{VB}$-groupoids have isomorphic cohomology. In Section 6 we study projectable sections of $\mathcal{VB}$-groupoids. This gives an alternative way to relate the complexes associated to Morita equivalent $\mathcal{VB}$-groupoids, and is needed in Section 7.

Similarly to the case of $\mathcal{VB}$-groupoids, we define Morita maps and a notion of Morita equivalence for $\mathcal{LA}$-groupoids.  We show in Theorem \ref{theorem: invariance of Lie 2-algebras} that if $\mathcal V$ and $\mathcal W$ are Morita equivalent $\mathcal{LA}$-groupoids then the corresponding Lie 2-algebras are isomorphic in the derived category of Lie 2-algebras. We note in Theorem \ref{theorem: L infinity map} that this implies that there exists an $L_\infty$ quasi-isomorphism between the two Lie 2-algebras, considered as 2-term differential graded Lie algebras.

As applications of our main results, we study the algebraic structure underlying multiplicative vector fields on a Lie groupoid. It turns out that the category of multiplicative vector fields on a Lie groupoid $\G$ is equipped with a Lie 2-algebra structure which dependes only on the Morita equivalence class of $\G$. This gives a positive answer to a conjecture by Hepworth \cite{hepworth}. Then, we use the Lie 2-algebra of multiplicative vector fields to introduce geometric vector fields on quotient stacks as degree one cohomology classes of multiplicative vector fields. This notion is more geometric rather than categorical and allows to describe vector fields on quotient stacks in many different situations, including classifying stacks, orbifolds, foliations among others.

The paper is organized as follows. In order to fix our notation, in Section 2 we briefly discuss the main examples of $\vb$-groupoids and the relation with representations up to homotopy. We also review the basics on $\la$-groupoids and the relation between Lie 2-algebras and crossed modules of Lie algebras. In Section 3 we introduce the category of multiplicative sections of a $\vb$-groupoid and we prove Theorem \ref{2vect}, which is the main result of this section. We also describe the category of multiplicative sections of an important class of $\vb$-groupoids, those which are regular. In Section 4 we study the cohomology of the complex of multiplicative sections of a $\vb$-groupoid and its relation with representations up to homotopy. In Section 5 we prove Morita invariance of the category of multiplicative sections of a $\vb$-groupoid, this is the content of Theorem \ref{moritainvarianceVB}. In Section 6 we study projectable sections of  $\mathcal{VB}$-groupoids and prove Theorem \ref{thm: fibre product of qisoms}, which is used in Section 7. In Section 7 we show that the category of multiplicative sections of an $\mathcal{LA}$-groupoid carries a natural Lie 2-algebra structure (Theorem \ref{Lie2}) and show that this structure is Morita invariant in an appropriate sense (Theorems \ref{theorem: invariance of Lie 2-algebras} and
\ref{theorem: L infinity map}). Finally, in Section \ref{sec:vectorfields}, we give a detailed exposition on the Lie 2-algebra of vector fields on a Lie groupoid. Several examples are discussed and we introduce a geometric notion of vector field on a differentiable stack.

While completing this work we learned that E. Lerman and D. Berwick-Evans have independently proved that the category of multiplicative vector fields on a Lie groupoid has a Lie 2-algebra structure, see \cite{BELerman}. In the case of tangent groupoids our results on $\mathcal{VB}$ and $\mathcal{LA}$-groupoids agree with theirs, but our methods are different.

\textbf{Acknowledgements:} The authors would like to thank A. Cabrera, O. Brahic, T. Drummond and M. Zambon for many useful comments and suggestions. Also, C.Ortiz thanks Newcastle University for the hospitality while part of this work was carried out. C.Ortiz was partially supported by grant 2016/01630-6 Sao Paulo Research Foundation and by Projeto Universal 482796/2013-8 CNPq-Brazil. J. Waldron thanks the University of Sao Paulo for their hospitality while part of this work was carried out. J. Waldron was supported by the EPSRC under the EPSRC Postdoctoral Prize Fellowship while part of this work was carried out.


\section{Preliminaries}

\subsection{2-vector spaces}
\label{subsection: 2-vector spaces}

2-vector spaces and their relation to 2-term complexes first appeaered in a more general setting in \cite{Deligne}. See  \cite{Daenzer} for an exposition of the (2-categorical) equivalence between 2-vector spaces and 2-term complexes, and \cite{BaezCrans} for details relevent to Lie 2-algebras. We denote by $\mathbf{Vect}$ the category of vector spaces over $\R$. A \textbf{2-vector space} is a category internal to $\mathbf{Vect}$. We depict a 2-vector space by a diagram of the form
\[
V_1 \rightrightarrows V_0
\]
where the horizontal arrows represent the source and target morphisms, $s$ and $t$.
A morphism of 2-vector spaces is an internal functor in $\mathbf{Vect}$. 2-vector spaces and morphisms between them form a category, denoted $\mathbf{2Vect}$. The forgetful functor $\mathbf{Vect} \to \mathbf{Set}$ extends to a faithfull functor $\mathbf{2Vect} \to \mathbf{Cat}$ allowing one to consider 2-vector spaces as small categories with some extra structure.  A morphism of 2-vector spaces is an {\bf equivalence} if the underlying functor is fully-faithful and essentially surjective. 

The category $\textbf{2-Term}$ is the category of cochain complexes of real vector spaces concentrated in degrees $0$ and $1$. We call the objects of $\textbf{2-Term}$ `2-term complexes'. There is an equivalence of categories:
\begin{equation}
\label{eqn: equivalence of 2Vect and 2Term}
\textbf{2-Vect} \simeq \textbf{2-Term}
\end{equation}
At the level of objects the equivalence (\ref{eqn: equivalence of 2Vect and 2Term}) is as follows. The functor $\textbf{2-Vect} \to \textbf{2-Term}$ is
\[
V_1 \rightrightarrows V_0  \; \;  \mapsto \; \;  \mathrm{Ker}(s) \xrightarrow{t|_{\mathrm{Ker}(s)}} V_0  
\]
and the quasi-inverse functor $\textbf{2-Term} \to \mathbf{2Vect}$ is
\[
V \xrightarrow{\partial} W  \; \; \mapsto \;  \; 
V \oplus W \rightrightarrows W
\]
where $V \oplus W \rightrightarrows W$ is the category associated to the action of $V$ on $W$ given by the map $(v,w) \mapsto \partial(v)+w$. Under the equivalence (\ref{eqn: equivalence of 2Vect and 2Term}), quasi-isomorphisms of 2-term complexes correspond to equivalences of 2-vector spaces, $H^0$ of a complex corresponds to the vector space of automorphisms of $0 \in V_0$, and $H^1$ corresponds to the vector space of isomorphism classes of objects in the category $V_1 \rightrightarrows V_0$.

Note that it is common, particularly in the related theory of 2-groups, that $\textbf{2-Term}$ is taken to be the category of chain complexes concentrated in degrees $1$ and $0$. Our convention is motivated by the link between our results and $\mathcal{VB}$-groupoid cohomology - see section \ref{sec:cohomology}.

\subsection{Lie 2-algebras and crossed modules of Lie algebras}
\label{subsection: Lie 2-algebras}

See \cite{BaezCrans} and \cite{Noohi} for the relevant aspects of the theory of Lie 2-algebras.
A {\bf Lie 2-algebra} is an internal category in $\mathbf{Lie}$. We depict a Lie 2-algebra by a diagram of Lie algebras
\[
\mathfrak g_1 \rightrightarrows \mathfrak g_0 
\]
where the horizontal arrows represent the source and target morphims. A morphism of Lie 2-algebras is an internal functor in $\mathbf{Lie}$. Lie 2-algebras and morphisms between them form a category, denoted $\mathbf{2Lie}$. The forgetful functor $\mathbf{Lie} \to \mathbf{Vect}$ extends to a functor $\mathbf{2Lie} \to \mathbf{2Vect}$ allowing one to consider Lie 2-algebras as 2-vector spaces with some extra structure.  Similarly, the forgetful functor $\mathbf{Lie} \to \mathbf{Set}$ extends to a functor $\mathbf{2Lie} \to \mathbf{Cat}$ allowing one to consider Lie 2-algebras as small categories with some extra structure. A morphism of Lie 2-algebras is an {\bf equivalence} if the underlying functor is fully-faithful and essentially surjective. Note that equivalences of Lie 2-algebras are in general not (quasi-) invertible.

A {\bf crossed module of Lie algebras}, or just {\bf crossed module}, consists of a pair of Lie algebras $\mathfrak g, \mathfrak h$ and Lie algebra morphisms
\begin{equation}
\label{eqn: X-mod}
\mathfrak g \xrightarrow{\partial} \mathfrak h \xrightarrow{\phi} \mathrm{Der}\left(\mathfrak g\right)
\end{equation}
such that
\begin{enumerate}
\item $\phi_{\partial (u)} = \mathrm{ad}_u$ for all $u \in \mathfrak g$.
\item $\partial \circ \phi_X   = \mathrm{ad}_X \circ \partial$ for all $X \in \mathfrak h$.
\end{enumerate}
We denote a crossed module by a diagram of type (\ref{eqn: X-mod}), or by a tuple $(\mathfrak g, \mathfrak h, \partial, \phi)$.  If $(\mathfrak g, \mathfrak h, \partial, \phi)$ is a crossed module then it follows from 1. that $H^0(\partial) = \mathrm{Ker}(\partial)\subseteq \mathfrak{g}$ is an abelian Lie subalgebra and from 2. that $\mathrm{Im}(\partial)$ is an ideal in $\mathfrak h$ and therefore that $H^1(\partial) = \mathrm{Coker}\left(\partial\right)$ has a natural Lie algebra structure.  Crossed modules form a category $\textbf{X-mod}$, where a morphism
\[
(\mathfrak g, \mathfrak h, \partial, \phi) 
\to
(\mathfrak g', \mathfrak h', \partial', \phi')
\]
consists of Lie algebra morphisms $f_1:\mathfrak g \to \mathfrak g'$, $f_2: \mathfrak h \to \mathfrak h'$, such that
\begin{itemize}
\item[i)] $f_2 \circ \partial = \partial' \circ f_1$.
\item[ii)] $f_1 \left( \phi_X(u) \right)  = \phi'_{f_2(X)} \left( f_1(u)\right)$ for all $X \in \mathfrak h, u \in \mathfrak g$.
\end{itemize}
\begin{remark} \label{remark: morphism of Xmods}
In fact, if $f_1$ is only assumed to be a linear map, it follows from the axioms of crossed modules that if i) and ii) hold then $f_1$ is automatically a Lie algebra morphism. Indeed, if $u,v \in \mathfrak g$ then
\begin{align*}
f_1 \left(   [u,v] \right)  & =  f_1 \left( \mathrm{ad}_u v \right) \\
& = f_1 \left(  \phi_{\partial(u)} (v) \right) \\
& = \phi'_{f_2 \partial(u)} \left( f_1(v) \right) \\
& = \phi'_{\partial' f_1(u)} \left( f_1(v)\right) \\
& = \mathrm{ad}_{f_1(u)} f_1(v) \\
& = \left[  f_1(u) , f_1(v) \right]
\end{align*}
\end{remark}
A morphism of crossed modules induces a pair of Lie algebra morphisms $\mathrm{Ker}\left(\partial \right) \to \mathrm{Ker}\left(\partial'\right)$ and $\mathrm{Coker}\left(\partial \right) \to \mathrm{Coker}\left(\partial'\right)$. A morphism of crossed modules is a {\bf quasi-isomorphism} if both of these Lie algebra morphisms are isomorphisms. There is a forgetful functor $\textbf{X-mod}\to \textbf{2-Term}$ defined on objects by
\[
\mathfrak g \xrightarrow{\partial} \mathfrak h \xrightarrow{\phi} \mathrm{Der}\left(\mathfrak g\right) 
\; \;  \mapsto \;  \;
\mathfrak g \xrightarrow{\partial} \mathfrak h 
\]
and which preserves the notion of quasi-isomorphism.

There is an equivalence of categories
\begin{equation}
\label{eqn: 2Lie and X-mod}
\mathbf{2Lie}  \simeq \textbf{X-mod}
\end{equation}
At the level of objects the equivalence is as follows. The functor $\mathbf{2Lie} \to \textbf{X-mod}$ is
\[
\mathfrak g_1 \rightrightarrows \mathfrak g_0  \; \;  \mapsto \; \;  \mathrm{Ker}(s) \xrightarrow{t|_{\mathrm{Ker}(s)}} \mathfrak g_0  \xrightarrow{\mathrm{ad} \circ 1} \mathrm{Der}\left( \mathrm{Ker}(s)\right)
\]
and the quasi-inverse functor $\textbf{X-mod} \to \mathbf{2Lie}$ is
\[
\mathfrak g \xrightarrow{\partial} \mathfrak h  \; \; \mapsto \;  \; \mathfrak g \rtimes \mathfrak h \rightrightarrows \mathfrak h
\]
where $\mathfrak g \rtimes \mathfrak h$ is the vector space $\mathfrak g \oplus \mathfrak h$ is equipped with the semi-direct product Lie bracket determined by $\phi$;
\[
\left[  \left(u,X\right) , \left(v,Y\right)  \right]_\phi  = 
\left(  [u,v] + \phi_X (v) - \phi_Y(u)  , [X,Y] \right)
\]
 and the structure morphisms of $\mathfrak g \rtimes \mathfrak h \rightrightarrows \mathfrak h$ are identical to those of the 2-vector space associated to the 2-term complex $\mathfrak g \xrightarrow{\partial} \mathfrak h$. Under the equivalence (\ref{eqn: 2Lie and X-mod}) a morphism in $\mathbf{2Lie}$ is an equivalence if and only if the corresponding morphism in $\textbf{X-mod}$ is a quasi-isomorphism. We shall occasionally use the equivalence (\ref{eqn: 2Lie and X-mod}) implicitly, and shall speak of Lie 2-algebras interchangeably. The \textbf{derived category of Lie 2-algebras} is the localisation of $\textbf{X-mod}$ (equivalently $\textbf{2Lie}$) obtained by inverting all quasi-isomorphisms (equivalently all equivalences).


We can summarise the discussion above in the following commutative diagram of categories and functors:

\[
\xymatrix{
\mathbf{2Lie} \ar@<+.5ex>[r]  \ar[d] &  \textbf{X-mod} \ar[d] \ar@<+.5ex>[l]\\
\mathbf{2Vect}  \ar@<+.5ex>[r] \ar[d]  &   \textbf{2-Term} \ar@<+.5ex>[l]  \\
\mathbf{Cat}   & 
}
\]
The vertical arrows are forgetfull functors, with those on the left preserving equivalences and the functor on the right preserving quasi-isomorphisms, and the horizontal functors are equivalences of categories, all of which map equivalences to quasi-isomorphisms or vice-versa. 

Note that the category $\textbf{X-mod}$ is also equivalent to the category of 2-term differential graded Lie algebras, though we will only make use of this equivalence briefly in Section \ref{subsection: Linfinity}.

\subsection{$\vb$-groupoids and representations up to homotopy}\label{subsec:vbgruths}

Let $\G\rr M$ be a Lie groupoid. A \textbf{$\vb$-groupoid} over $\G$ is a Lie groupoid $\mathcal{V}\rr E$, where both $q_{\mathcal{V}}:\mathcal{V}\rmap \G$ and $q_E:E\rmap M$ are vector bundles which are compatible with the groupoid structures in the sense that all the structure maps of the Lie groupoid $\mathcal{V}\rr E$ are vector bundle morphism over the corresponding structure maps of $\G\rr M$. The structure maps of $\mathcal{V}\rr E$ are denoted by $\tilde{s},\tilde{t},\tilde{m},\tilde{i}$ and $\tilde{1}$. We write a $\vb$-groupoid as a square

\begin{align}\label{VBgroupoid}
\SelectTips{cm}{}\xymatrix{ \mathcal{V} \ar@<0.5ex>[r]^{\tilde{s}}\ar@<-.5ex>[r]_{\tilde{t}}\ar[d]_{q_{\mathcal{V}}} &E \ar[d]^{q_E}\\
           \G \ar@<0.5ex>[r]^{s}\ar@<-.5ex>[r]_{t} & M}
\end{align}

The fact that $\tilde{s}:\mathcal{V}\rmap E$ is a surjective submersion implies that $\ker(\tilde{s})\subseteq \mathcal{V}$ is a vector bundle over $\G$. The \textbf{core} of $\mathcal{V}$ is the vector bundle over $M$ defined as the pullback $C:=\ker(\tilde{s})|_{M}$. The  core $C$ fits into an exact sequence of vector bundles over $\G$:

\begin{equation}\label{coresequence}
\xymatrix{0\ar[r]&t^*C\ar[r]^{r}&\mathcal{V}\ar[r]^{\tilde{s}}&s^*E\ar[r]&0},
\end{equation}

\noindent where the map $r:t^*C\rmap \mathcal{V}$ is defined by $r(c_{t(g)})=c_{t(g)}\tilde{0}_g$. We refer to \eqref{coresequence} as the \textbf{core sequence} of $\mathcal{V}$. Sections of the core determine special sections of $\mathcal{V}$, as we explained below.

\begin{definition}\label{rightinvariant}

Let $c\in \Gamma(C)$ be a section of the core. The \textbf{right invariant} section $c^r\in\Gamma(\mathcal{V})$ determined by  $c$ is defined as $c^r(g)=c(t(g))\tilde{0}_g$ for $g\in\G$. Similarly, the \textbf{left invariant} section $c^l\in\Gamma(\mathcal{V})$ is defined by $c^l(g)=-\tilde{0}_g\tilde{i}(c(s(g)))$, $g\in\G$.

\end{definition}

\begin{example}(2-vector spaces)\label{2vectsVB}
A $\vb$-groupoid over the trivial Lie groupoid $\{*\}\rr \{*\}$ is the same as a finite dimensional 2-vector space.
\end{example}

\begin{example}(Tangent groupoid)\label{tangentVB}
Let $\G\rr M$ a Lie groupoid. The \textbf{tangent groupoid} is the Lie groupoid $T\G\rr TM$ obtained by applying the tangent functor to each of the structure maps of $\G$. The tangent groupoid is a $\vb$-groupoid with core $A$ the Lie algebroid of $\G$.  
\end{example}

\begin{example}(Cotangent groupoid)\label{cotangentVB}
Given a Lie groupoid $\G\rr M$ with Lie algebroid $A$, the cotangent bundle $T^*\G$ inherits a Lie groupoid structure over the dual $A^*$, making $T^*\G\rr A^*$ into a $\vb$-groupoid with core $T^*M$. This is the \textbf{cotangent groupoid} of $\G$. For more details see \cite{CDW}.
\end{example}

\begin{example}(Representations)\label{transformationVB}
A \textbf{representation} of a Lie groupoid $\G\rr M$ is a vector bundle $E\rmap M$ equipped with a linear action of $\G$, i.e. a morphism of Lie groupoids $\Delta^E:\G\rmap GL(E)$, where $GL(E)\rr M$ is the Lie groupoid with objects being points of $M$ and arrows $x\rmap y$ defined as linear isomorphisms between the fibers $E_x\rmap E_y$. Any representation gives rise to a $\vb$-groupoid by considering the transformation groupoid $s^*E=\G\ltimes E\rr E$. This is a $\vb$-groupoid with core $C=0$. Indeed, every $\vb$-groupoid with zero core arises in this way. 
\end{example}

\begin{example}\label{equivariantVB}(Equivariant vector bundles)
Let $E\rmap M$ be a vector bundle. Assume that $G$ is a Lie group which acts on $E$ by vector bundle automorphisms $\phi_g:E\rmap E$ covering a diffeomorphism $\psi_g:M\rmap M$. In particular, $G$ acts on $M$ via $\psi_g:M\rmap M$. The transformation groupoid $G\ltimes E\rr E$ has a natural structure of $\vb$-groupoid over $G\ltimes M\rr M$. Indeed, a representation of $G\ltimes M\rr M$ is equivalent to a $G$-equivariant vector bundle $E\rmap M$. Hence, the construction of Example \ref{transformationVB} applies and one gets a $\vb$-groupoid $s^*E\rr E$, where $s:G\ltimes M\rr M$ is the source map. One easily observes that the groupoid $s^*E\rr E$ coincides with $G\ltimes E\rr E$.
\end{example}

The notion of representation of a Lie groupoid can be generalized by representing Lie groupoids on graded vector bundles. This corresponds to the concept of representation up to homotopy introduced in \cite{AriasCrainic2}. For that, recall that a \textbf{quasi-action} of a Lie groupoid $\G\rr M$ on a vector bundle $E\rmap M$ is a smooth map $\Delta:\G\to L(E)$ with $\Delta_g:E_{s(g)}\rmap E_{t(g)}$ for any $g\in \G$. Here $L(E)$ is the smooth category whose objects are elements of $M$, and morphisms between $x,y\in M$ are linear maps $E_x\to E_y$ between the fibers. If $\Delta:\G\rmap L(E)$ satisfies $\Delta_x=Id:E_x\rmap E_x$, then we say that $\Delta$ is \textbf{unital}. Also, if $\Delta_{gh}=\Delta_g\circ \Delta_h$ for any composable elements $g,h\in \G$, then $\Delta$ is called \textbf{flat}. Hence, a representation of $\G$ is just a flat unital quasi-action.

A representation up to homotopy of a Lie groupoid $\G$ on a graded vector bundle $\mathcal{E}$ can be defined in cohomological terms, as a differential in the complex $C(G,\mathcal{E})$ of groupoid cochains with values in $\mathcal{E}$. In this work we are concerned only with 2-term representations up to homotopy, i.e. on a graded bundle of the form $C\oplus E$. For convenience, we use the following equivalent definition of representation up to homotopy.

\begin{definition}\label{def:ruthgroupoids}

Let $\G\rr M$ be a Lie groupoid. A \textbf{representation up to homotopy} of $G$ on the graded vector bundle $\mathcal{E}=C\oplus E$, is given by a quadruple $(\partial,\Delta^C,\Delta^E,\Omega)$ where:
\begin{itemize}
\item $\partial:C\to E$ is a bundle map, 
\item $\Delta^C$ and $\Delta^E$ are unital quasi-actions on $C$ and $E$, respectively, 
\item $\Omega\in \Gamma (\G^{(2)}, \Hom(s^*E,t^*C))$ is a section assigning to each composable pair $(g_1,g_2)\in \G^{(2)}$ a linear map $\Omega_{g_1,g_2}:E_{s(g_2)} \to C_{t(g_1)}$ which is normalized, i.e. $\Omega_{g_1,g_2}=0$ if either $g_1$ or $g_2$ is a unit, 
\end{itemize}
satisfying the following conditions:

\begin{align}
\Delta^E_{g_1}\circ \partial&=\partial\circ \Delta^C_{g_1}\nonumber\\
\Delta^C_{g_1}\Delta^C_{g_2}- \Delta^C_{g_1g_2} &+ \Omega_{g_1,g_2}\partial=0\nonumber\\
\Delta^E_{g_1}\Delta^E_{g_2}- \Delta^E_{g_1g_2} &+ \partial\Omega_{g_1,g_2}=0\nonumber\\
\Delta^C_{g_1}\Omega_{g_2,g_3}-\Omega_{g_1g_2,g_3}&+\Omega_{g_1,g_2g_3}-\Omega_{g_1,g_2}\Delta^E_{g_3}=0 \label{eq:repGeqs}
\end{align}
for every composable triple $(g_1,g_2,g_3)$ in $G$.

\end{definition}

As explained in \cite{GM10}, a representation up to homotopy of $\G$ on $C\oplus E$, induces a groupoid structure on the vector bundle $t^*C\oplus s^*E\rmap \G$ which makes $t^*C\oplus s^*E$ into a $\vb$-groupoid over $\G$ with core $C$ and side $E$. The groupoid $t^*C\oplus s^*E\rr E$ is defined as follows:

\begin{align}
 \tilde{s}(c,g,e)=e, \quad \tilde{t}(c,g,e)=\partial(c)+\Delta^E_g(e), \quad \tilde{1}_e=\left(0^C(x),1(x), e\right), \ e\in E_x \nonumber \\
(c_1,g_1,e_1)\cdot (c_2,g_2,e_2)=\left(c_1+\Delta^C_{g_1}(c_2)-\Omega_{g_1,g_2}(e_2),g_1\cdot g_2,e_2\right)
\nonumber \\
(c,g,e)^{-1}=(-\Delta^C_{g^{-1}}(c)+\Omega_{g^{-1},g}(e),g^{-1},\partial(c)+\Delta^E_g(e)). 
\label{eq:splitvbgmaps}
\end{align}

A $\vb$-groupoid of the form $t^*C\oplus s^*E\rr E$ with structure maps \eqref{eq:splitvbgmaps} is called a \textbf{split} $\vb$-groupoid.

Conversely, given a $\vb$-groupoid, the \textbf{core anchor map} is the vector bundle morphism $\partial:C\rmap E$ defined by $\partial(c)=\tilde{t}(c)$. In order to obtain a representation up to homotopy out of a $\vb$-groupoid, first we need to choose a horizontal lift. Note that, for every unit element $x\in M$, the fiber $\mathcal{V}_x$ is canonically isomorphic to $C_x\oplus E_x$ via the isomorphism 

\begin{equation}\label{canonicalunits}
C_x\oplus E_x\rmap \mathcal{V}_x; (c,e)\mapsto c+\tilde{1}(e).
\end{equation}

\noindent As defined in \cite{GM10}, a \textbf{horizontal lift} of $\mathcal{V}$ is a right splitting $h:s^*E\rmap \mathcal{V}$ of \eqref{coresequence} which at every unit $x\in M$ coincides with the canonical splitting \eqref{canonicalunits}. As shown in \cite{GM10}, horizontal lifts of a $\vb$-groupoid always exist. The choice of a horizontal lift on $\mathcal{V}$ induces a representation up to homotopy of $\G$ on $C\oplus E$. Actually, a result of Gracia Saz and Mehta \cite{GM10} says that there is a one-to-one correspondence between representations up to homotopy of $\G$ on $C\oplus E$ and $\vb$-groupoids $\mathcal{V}$ with core $C$ and side $E$ together with a horizontal lift. More precisely, as shown in \cite{dHO}, this construction is functorial and yields an equivalence of categories between the category of $\vb$-groupoids over $\G\rr M$ and that of representations up to homotopy of $\G\rr M$.

\begin{remark}
In the special case of the Lie groupoid $\{*\}\rr \{*\}$ a representation up to homotopy of $\{*\}\rr \{*\}$ is just a linear map between vector spaces $\partial: V_0\rmap V_1$ and one recovers the (finite dimensonal version of the) equivalence of categories $\mathbf{2Vect}\simeq \mathbf{2Term}$ described in Section \ref{subsection: 2-vector spaces} above.
\end{remark}

It will be useful to have an expression for the right and left invariant sections associated to $c\in\Gamma(C)$ in the case of a split $\vb$-groupoid. This follows straightforward from the definition of the structure maps \eqref{eq:splitvbgmaps}, yielding

\begin{equation}\label{splitrightinvariant}
c^r(g)=c_{t(g)}\cdot \tilde{0}_g=(c_{t(g)},g,0),
\end{equation}

and

\begin{equation}\label{splitleftinvariant}
c^l(g)=-\tilde{0}_g\cdot \tilde{i}(c_{s(g)})=(\Delta^C_gc_{s(g)},g,-\partial c_{s(g)}).
\end{equation}

\subsection{$\la$-groupoids}\label{LAgroupoids}

In this subsection we follow \cite{mackII} closely. Let $\G\rr M$ be a Lie groupoid. An \textbf{$\la$-groupoid} over $\G$ is a $\vb$-groupoid as in \eqref{VBgroupoid} where both vertical structures $\mathcal{V}\rmap \G$ and $E\rmap M$ are Lie algebroids and all the structure maps defining the Lie groupoid $\mathcal{V}\rr E$ are Lie algebroid morphisms over the corresponding structure maps of $\G\rr M$.

The core $C\rmap M$ of an $\la$-groupoid inherits a Lie algebroid structure defined in the following way. The anchor map $\rho_C:C\rmap TM$ is given by $\rho_C(c)=\rho_E\circ \partial$, where $\rho_E:E\rmap TM$ is the anchor of $E$ and $\partial:C\rmap E$ is the core anchor map. Given two sections $c_1,c_2\in\Gamma(C)$ the Lie bracket $[c_1,c_2]_C\in\Gamma(C)$ is defined as the unique section of $C$ satisfying $[c_1,c_2]^r_C=[c^r_1,c^r_2]_{\mathcal{V}}$.

We finish this section with some examples of $\la$-groupoids.

\begin{example}\label{ex:lie2}

An $\la$-groupoid $V\rr W$ over $\{*\}\rr \{*\}$ is the same as a finite dimensional Lie 2-algebra.

\end{example}

\begin{example}

The tangent groupoid $T\G\rr TM$ is an $\la$-groupoid with respect to the canonical Lie algebroid structures on the tangent bundles $T\G\rmap \G$ and $TM\rmap M$. 

\end{example}

\begin{example}

A \textbf{Poisson groupoid} \cite{Wco} is a Lie groupoid $\G\rr M$ equipped with a Poisson structure $\pi\in\mathfrak{X}^2(\G)$ which is multiplicative, i.e. the graph of the multiplication $\Gamma_m=\{(g,h,gh);s(g)=t(h)\}\subseteq \G\times\G\times\bar{\G}$ is a coisotropic submanifold. As shown in \cite{MX1}, the multiplicativity condition is equivalent to $\pi^{\sharp}:T^*\G\rmap T\G$ being a Lie groupoid morphism. In this case, the cotangent groupoid $T^*\G\rr A^*$ inherits the structure of an $\la$-groupoid. See \cite{mackII} for more details. More generally, multiplicative Dirac structures \cite{Ortiz13} can be viewed as $\la$-groupoids $L\rr E$ with respect to the natural structure induced from the Courant algebroid $T\G\oplus T^*\G$ which is a Lie groupoid over $TM\oplus A^*$.
\end{example}

\begin{example}
Suppose that $\Gamma$ is a discrete group acting on a Lie algebroid $E \to M$ by Lie algebroid automorphisms. We observed in Example ... that the transformation groupoid $\Gamma \ltimes E \rightrightarrows E$ is a $\mathcal{VB}$-groupoid over $\Gamma \ltimes M \rightrightarrows M$. The vector bundle $\Gamma \ltimes E = s^\ast E$ has a natural Lie algebroid structure induced by that of $E$, with respect to which $\Gamma \ltimes E \rightrightarrows E$ is an $\mathcal{LA}$-groupoid over $\Gamma\ltimes M \rightrightarrows M$.
\end{example}



\section{The category of multiplicative sections}\label{categorymultiplicative}

In this section we introduce the category of multiplicative sections of a $\vb$-groupoid. We study the algebraic structure underlying such category in both cases of $\vb$-groupoids and $\la$-groupoids.

\subsection{Multiplicative sections}

Let $\mathcal{V}\rr E$ be a $\vb$-groupoid over $\G\rr M$. A \textbf{multiplicative section} of $\mathcal{V}$ is given by a pair of sections $e\in\Gamma(E)$ and $V\in\Gamma(\mathcal{V})$ such that $V:\G\rmap \mathcal{V}$ is a Lie groupoid morphism covering $e:M\rmap E$. The space of multiplicative sections of $\mathcal{V}$ is denoted by $\Gamma_{mult}(\mathcal{V})$. Multiplicative sections are closed under addition and scalar multiplication, making $\Gamma_{mult}(\mathcal{V})$ into a subspace of $\Gamma(\mathcal V)$.

\begin{example}
If $V_1 \rightrightarrows V_0$ is a 2-vector space, considered as a $\mathcal{VB}$-groupoid over $\left\{\ast\right\} \rightrightarrows \left\{\ast\right\}$, then $\Gamma_{mult}\left(V_1 \rightrightarrows V_0\right)$ can be identified with $V_0$. Indeed, functors from the trivial groupoid $\left\{\ast\right\} \rightrightarrows \left\{\ast\right\}$ to $V_1 \rightrightarrows V_0$ can be identified with the set of units $V_0$ in $V_1 \rightrightarrows V_0$.
\end{example}

\begin{example}
In the case of the tangent groupoid $T\G\rr TM$, a multiplicative section is referred to as a \textbf{multiplicative vector field} on $\G$. The space of multiplicative vector fields on $\G$ is denoted by $\mathfrak{X}_{mult}(\G)$. 
\end{example}

\begin{example}
Multiplicative sections of the cotangent groupoid $T^*\G\rr A^*$ are called \textbf{multiplicative 1-forms} on $\G$. The space of multiplicative 1-forms on $\G$ is denoted by $\Omega^1_{mult}(\G)$. This is in agreement with the terminology introduced in \cite{MX3}.
\end{example}

\begin{example}
\label{example: sections of transformationVB}
If $\mathcal G \ltimes E \rightrightarrows E$ is a $\mathcal{VB}$-groupoid associated to a representation of a Lie groupoid $\mathcal G \rightrightarrows M$ on a vector bundle $E \to M$ (see Example \ref{transformationVB}) then $\Gamma_{mult}\left(\mathcal G \ltimes E\right)$  can be identified with the vector space $\Gamma\left(E\right)^\mathcal G$ of $\mathcal G$-invariant sections of $E$:
\[
\Gamma_{mult}\left(\mathcal G \ltimes E\right) \cong 
\Gamma_{mult}\left(E\right)^\mathcal G
= 
\left\{  e \in \Gamma(E) \; | \;  g \cdot e(s(g)) = e(t(g)) \; \forall \; g \in \mathcal G \right\}
\]
Indeed, if $\left(V,e\right)$ is a multiplicative section of $\mathcal G \ltimes E$ then $\tilde s \circ V = e \circ s$ implies that $V=s^\ast e$ and $\tilde t \circ V = e \circ t$ then implies that $e \in \Gamma\left(E\right)^\mathcal G$. An easy calculation shows the converse - that if $e \in \Gamma\left(E\right)^\mathcal G$ then the pair $\left(s^\ast e,e\right)$ is automatically multiplicative.
\end{example}

\begin{example}
As a special case of the previous example, if $G \ltimes E \rightrightarrows E$ is the $\mathcal{VB}$-groupoid associated to a $G$-equivariant vector bundle $E\to M$ as in Example \ref{equivariantVB} then $\Gamma\left(E\right)^{G \ltimes E} = \Gamma\left(E\right)^G$ and therefore
\[
\Gamma_{mult}\left(G \ltimes E\right) \cong \Gamma\left(E\right)^G
\]
\end{example}

One easily observes that a $\vb$-groupoid morphism $\Phi:\mathcal{V}_1\rmap \mathcal{V}_2$ covering $id:\G\rmap \G$, induces a map $\Gamma(\mathcal{V}_1)\rmap \Gamma(\mathcal{V}_2)$, which preserves multiplicative sections. Moreover, if $\Phi$ is a $\vb$-groupoid isomorphism, then $\Phi$ induces an isomorphism of vector spaces $\Gamma_{mult}(\mathcal{V}_1)\cong \Gamma_{mult}(\mathcal{V}_2)$. In particular, given a $\vb$-groupoid $\mathcal{V}$ as in \eqref{VBgroupoid} together with a horizontal lift $h:s^*E\rmap \mathcal{V}$, one gets a representation up to homotopy of $\G$ on $C\oplus E$. This determines a $\vb$-groupoid isomorphism between $\mathcal{V}$ and a split $\vb$-groupoid $t^*C\oplus s^*E$ described in \eqref{eq:splitvbgmaps}. As a consequence, multiplicative sections of $\mathcal{V}$ are in one-to-one correspondence with multiplicative sections of a split $\vb$-groupoid.

\begin{proposition}

Let $\mathcal{V}$ be a $\vb$-groupoid over $\G\rr M$ with core $C$. For any section $c\in \Gamma(C)$, the section $c^r-c^l:\G\rmap \mathcal{V}$ is multiplicative.

\end{proposition}

\begin{proof}
An isomorphism of $\vb$-groupoids $\mathcal{V}\rmap \mathcal{V}'$ covering $id:\G\rmap \G$ sends left and right invariant sections to left and right invariant sections respectively. Hence, we may assume that $\mathcal{V}=t^*C\oplus s^*E$ is a split $\vb$-groupoid associated to a representation up to homotopy $(\partial,\Delta^C,\Delta^E,\Omega)$ of $\G$ on $C\oplus E$. The structure maps of $t^*C\oplus s^*E\rr E$ are given by \eqref{eq:splitvbgmaps}. Consider now a section $c\in \Gamma(C)$, then equations \eqref{splitrightinvariant} and \eqref{splitleftinvariant} describe the corresponding right and left invariant sections associated to $c$. Hence,

\begin{align*}
(c^r-c^l)(g)=&(c_{t(g)},g,0)-(\Delta^C_gc_{s(g)},g,-\partial c_{s(g)})\\
=&(c_{t(g)}-\Delta^C_gc_{s(g)},g,\partial c_{s(g)}).
\end{align*}

On the one hand, let $g,h\in \G$ be composable arrows, i.e. $s(g)=t(h)$. Then,

\begin{align*}
(c^r-c^l)(gh)=&(c_{t(g)}-\Delta^C_{gh}c_{s(h)},gh,\partial c_{s(h)})\\
=& (c_{t(g)}-\Delta^C_g\Delta^C_hc_{s(h)}-\Omega_{g,h}\partial c_{s(h)},gh,\partial c_{s(h)}),
\end{align*}

\noindent where the last equality follows from the second condition of a representation up to homotopy \eqref{eq:repGeqs}. On the other hand, we have that

\begin{align*}
(c^r-c^l)(g)\cdot (c^r-c^l)(h)=& (c_{t(g)}-\Delta^C_gc_{s(g)},g,\partial c_{s(g)})\cdot (c_{t(h)}-\Delta^C_hc_{s(h)},h,\partial c_{s(h)})\\
=& (c_{t(g)}-\Delta^C_gc_{s(g)}+\Delta^C_gc_{t(h)}-\Delta^C_g\Delta^C_hc_{s(h)}-\Omega_{g,h}\partial c_{s(h)},gh,\partial c_{s(h)})\\
=&(c_{t(g)}-\Delta^C_g\Delta^C_hc_{s(h)}-\Omega_{g,h}\partial c_{s(h)},gh,\partial c_{s(h)}),
\end{align*}

\noindent where in the previous identities we have used the second formula in \eqref{eq:splitvbgmaps} for the multiplication on a split $\vb$-groupoid, together with the identity $\Delta^C_gc_{s(g)}=\Delta^C_gc_{t(h)}$, which holds because $s(g)=t(h)$. Hence, we conclude that for any composable arrows $g,h\in \G$, the following identity holds

$$(c^r-c^l)(gh)=(c^r-c^l)(g)\cdot (c^r-c^l)(h).$$

\noindent This proves that $c^r-c^l$ is a multiplicative section.
\end{proof}

\subsection{The category of multiplicative sections}

Consider a $\vb$-groupoid $\mathcal{V}\rr E$ over $\G\rr M$. We will introduce the category of multiplicative sections of $\mathcal{V}$. Since multiplicative sections are special type of functors $\G\rmap \mathcal{V}$, one can talk about (smooth) natural transformations between multiplicative sections. We look at those natural transformations which are also compatible with the linear structure. On any $\vb$-groupoid \eqref{VBgroupoid}, the bundle projection $q_{\mathcal{V}}:\mathcal{V}\rmap \G$ is a groupoid morphism over $q_E:E\rmap M$ (see \cite{GM10}). Given multiplicative sections $V,V':\G\rmap \mathcal{V}$ a \textbf{morphism} is a natural transformation $\tau:V\Rightarrow V'$ such that $1_{q_{\mathcal{V}}} \bullet \tau=1_{id_{\G}}$, where $\bullet$ denotes the horizontal composition of natural transformations:
\[\xymatrix@!@C=50pt{\G \ar@/_1.1pc/[r]|{\overset{ }{\, V\, }}="a"  \ar@/^1.1pc/[r]|{\, {V'}_{\,} \, }="b"
                  & \mathcal{V} \ar@/_1.1pc/[r]|{\overset{ }{\, q_{\mathcal{V}}\, }}="c"  \ar@/^1.1pc/[r]|{\, {q_{\mathcal{V}}} \, }="d" 
							  	  & \G \ar@{=>}^-{\tau}"a";"b"-<0pt,4pt>\ar@{=>}^-{1_{q_{\mathcal{V}}}}"c";"d"-<0pt,4pt>}
\quad = \quad
\xymatrix@!@C=50pt{\G \ar@/_1.1pc/[r]|{{ }{\, id_{\mathcal G}\, }}="a"  \ar@/^1.1pc/[r]|{\, id_{\mathcal G} \, }="c" & \G \ar@{=>}^{1_{id_{\mathcal G}}}"a";"c"-<0pt,4pt>}\]

\begin{definition}
Let $\mathcal{V}$ be a $\vb$-groupoid as in \eqref{VBgroupoid}. The \textbf{category of multiplicative sections} of $\mathcal{V}$ is defined as the subcategory $\mathrm{Sec}(\G,\mathcal{V})$ of $\Hom_{\mathrm{LieGpd}} (\G, \mathcal{V})$ defined by:
\begin{itemize}
\item[i)] the objects of $\mathrm{Sec}(\G,\mathcal{V})$ are multiplicative sections $V: \G\rmap \mathcal{V}$,

\item[ii)] the morphisms $\tau : V \Rightarrow W$ of $\mathrm{Sec}(\G,\mathcal{V})$ are natural transformations $\tau$ such that $1_{q_{\mathcal{V}}} \bullet \tau=1_{id_{\G}}$. 
\end{itemize}
\end{definition}

\begin{proposition} \label{proposition: description of morphisms between sections}
Let $\mathcal{V}$ be a $\vb$-groupoid as in \eqref{VBgroupoid}. Let $V,V' \in \mathrm{Sec}(\G,\mathcal{V})$ be multiplicative sections and $\tau: V \Rightarrow V'$ an arbitrary natural transformation. Then $\tau$ is a morphism in $\mathrm{Sec}(\G,\mathcal{V})$ if and only if $\tau(x) \in \mathcal{V}_{1_x}$ for all $x \in M$.
\end{proposition}

\begin{proof}
Let $x \in M$ a unit of $\G$. We have
\begin{align*}
\left( 1_{q_{\mathcal{V}}} \bullet \tau \right) (x) & = 1_{q_{\mathcal{V}}} \left( q_{E} \circ e'(x) \right) \circ q_{\mathcal{V}} \left( \tau(x) \right)  \\
& = 1_{q_{\mathcal{V}}} (x) \circ q_{\mathcal{V}} \left (\tau(x) \right) \\
& = 1_x \circ q_{\mathcal{V}} \left( \tau(x) \right) \\
& = q_{\mathcal{V}} \left( \tau(x) \right)  
\end{align*}
whereas
\begin{align*}
1_{id_{\G}} (x) & = 1_x
\end{align*}
Therefore $1_{q_{\mathcal{V}}} \bullet \tau=1_{id_{\G}}$ if and only if $q_{\mathcal{V}} \circ \tau = 1:M\rmap \G$, or equivalently, if and only if $\tau(x) \in \mathcal{V}_{1_x}$ for all $x \in M$.
\end{proof}

In other words, the objects of $\mathrm{Sec}(\G,\mathcal{V})$ are multiplicative sections, and by Proposition \ref{proposition: description of morphisms between sections}, we can identify morphisms in $\mathrm{Sec}(\G,\mathcal{V})$ with certain sections of the vector bundle $1^* \mathcal{V} \cong C \oplus E$.

\subsection{Description of morphisms}

Let $\mathcal{V}$ be a $\vb$-groupoid as in \eqref{VBgroupoid}. Consider the category $\mathrm{Sec}(\G,\mathcal{V})$ of multiplicative sections of $\mathcal{V}$. By Proposition \ref{proposition: description of morphisms between sections}, a morphism $\tau:V\Rightarrow V'$ in $\mathrm{Sec}(\G,\mathcal{V})$ can be identified with a special section of $1^* \mathcal{V} \cong C \oplus E$, so we can write $\tau=(c_0,e_0)$ where $c_0\in\Gamma(C)$ and $e_0\in \Gamma(E)$. The following result identifies precisely which sections $\tau$ arise this way.

\begin{proposition}\label{explicitmorphism}

Let $V, V'\in \Gamma_{mult}(\mathcal{V})$ be multiplicative sections covering $e,e'\in\Gamma(E)$, respectively.  Suppose that $\tau:M\rmap \mathcal{V}$ is a smooth map with $q_{\mathcal{V}}\circ \tau=1:M\rmap \G$. Identify $\tau$ with a section $(c_0,e_0)\in \Gamma(1^*\mathcal{V})\cong \Gamma(C)\oplus \Gamma(E)$. Then $\tau$ defines a morphism $\tau:V\Rightarrow V'$ in $\mathrm{Sec}(\G,\mathcal{V})$ if and only if the following conditions are fulfilled

\begin{enumerate}

\item $e_0=e,$

\item $e'=e_0+\partial(c_0),$

\item $V'=V+c^r_0-c^l_0$

\end{enumerate}

\end{proposition}

\begin{remark} Consider $\vb$-groupoids $\mathcal{V}_1$ and $\mathcal{V}_2$ over $\G$. An isomorphism of $\vb$-groupoids $\Phi:\mathcal{V}_1\rmap \mathcal{V}_2$ covering $id:\G\rmap \G$ induces an isomorphism of categories $\Gamma(\G,\mathcal{V}_1)\rmap \Gamma(\mathcal{V}_2)$. Therefore, in order to describe $\mathrm{Sec}(\G,\mathcal{V})$ in a more explicit way, we may assume that $\mathcal{V}$ is a split $\vb$-groupoid $\mathcal{V}=t^*C\oplus s^*E\rr E$ associated to a representation up to homotopy $(\partial,\Delta^C,\Delta^E,\Omega)$ of $\G$ on $C\oplus E$. Recall that the structure maps of $t^*C\oplus s^*E\rr E$ are defined by \eqref{eq:splitvbgmaps}. This is going to be used along the proof of Proposition \ref{explicitmorphism}.
\end{remark}

\begin{proof}

Let $t^*C\oplus s^*E\rr E$ be the $\vb$-groupoid associated to a representation up to homotopy $(\partial,\Delta^C,\Delta^E,\Omega)$ of $\G$ on $C\oplus E$. 
A smooth map $\tau:M\rmap \mathcal{V}$ is a natural transformation $V\Rightarrow V'$ if and only if the following conditions hold:

\begin{itemize}

\item[i)] $\tilde{s}\circ \tau=e,$
\item[ii)] $\tilde{t}\circ \tau=e',$
\item[iii)] $\tau(g)$ is natural in $g$ with respect to $V$ and $V'$.

\end{itemize}

For any $x\in M$ we have

\begin{equation*}
\tilde{s}(\tau(x))=\tilde{s}(c_0(x),x,e_0(x))=e_0(x), 
\end{equation*}

and

\begin{equation*}
\tilde{t}(\tau(x))=\tilde{t}(c_0(x),x,e_0(x))=\partial(c_0(x))+e_0(x).
\end{equation*}

\noindent Hence, condition 1. and 2. in the proposition are equivalent to i) and ii), respectively.

Now, assume that $\tau:M\rmap \mathcal{V}$ satisfies i) and ii). Given $g\in \G$, the naturality square for $\tau$ with respect to $V$ and $V'$ is

\begin{equation}\label{naturality}
\xymatrix{
e_{s(g)}   \ar[d]_{\tau(s(g))}  \ar[r]^{V(g)}  &   e_{t(g)}   \ar[d]^{\tau(t(g))}  \\
e'_{s(g)}   \ar[r]_{V'(g)}  &   e'_{t(g)}
}
\end{equation}

\noindent and $\tau(t(g))\cdot V(g)=V'(g)\cdot \tau(s(g))$.

Since $V:\G\rmap t^*C\oplus s^*E$ is a multiplicative section covering $e\in\Gamma(E)$, then at each point $g\in\G$ we write $V(g)=(c_{t(g)},g,e_{s(g)})$. Similary, we write $V'(g)=(c'_{t(g)},g,e'_{s(g)})$. 

On the one hand

\begin{align*}
\tau(t(g))V(g)&=(c_0(t(g)),t(g),e_0(t(g)))\cdot(c_{t(g)},g,e_{s(g)})\\
&= (c_0(t(g))+c_{t(g)},g,e_{s(g)})\\
&= (c_0(t(g)),g,0)+(c_{t(g)},g,e_{s(g)})\\
&= c^r_0(g)+V(g),
\end{align*}

\noindent where in the last equality we have used \eqref{splitrightinvariant}.

\noindent On the other hand,

\begin{align*}
V'(g)\cdot \tau(s(g))&=(c'_{t(g)},g,e'_{t(g)})\cdot (c_0(s(g)),s(g),e_0(s(g)))\\
&= (c'_{t(g)}+\Delta^C_gc_0(s(g)),g, e_0(s(g)))\\
&= (c'_{t(g)}+\Delta^C_gc_0(s(g)),g, e'_{s(g)}-\partial c_0(s(g)))\\
&= (c'_{t(g)},g,e'_{s(g)})+(\Delta^C_gc_0(s(g)),g,-\partial c_0(s(g)))\\
&= V'(g)+c^l_0(g),
\end{align*}

\noindent where in the last equality we used \eqref{splitleftinvariant}. Now, the commutativity of \eqref{naturality} is equivalent to $V'(g)=V(g)+c^r_0(g)-c^l_0(g)$, proving condition 3.
\end{proof}

\begin{definition}\label{complexmultiplicativesections}
Let $\mathcal{V}\rr E$ be a $\vb$-groupoid over $\G\rr M$ with core $C$. The \textbf{complex of multiplicative sections} of $\mathcal{V}$, denoted $C^\bullet_{mult}\left(\mathcal V\right)$, is the 2-term complex of vector spaces
\[
\delta: \Gamma(C)\to \Gamma_{mult}(\mathcal{V})
\]
 where 
 \[
 \delta(c)=c^r-c^l
 \]
\end{definition}

As recalled in Section \ref{subsection: 2-vector spaces}, associated to the 2-term complex $\delta: \Gamma(C) \to \Gamma_{mult}(\mathcal V)$ is a 2-vector space
\[
\Gamma(C) \oplus \Gamma_{mult}(\mathcal V) \rightrightarrows 
\Gamma_{mult}(\mathcal V)
\]

\begin{theorem}\label{2vect}
The category of multiplicative sections $\mathrm{Sec}(\G,\mathcal{V})$ is isomorphic to (the category underlying) the 2-vector space $\Gamma(C)\oplus \Gamma_{mult}(\mathcal{V})\rr \Gamma_{mult}(\mathcal{V})$, whose structure maps are
\[
\hat{s}(c,V)=V; \quad \hat{t}(c,V)=V+c^r-c^l
\]
\end{theorem}

\begin{proof}

We have already described a map $\mathrm{Sec}(\G,\mathcal{V})\rmap \Gamma(C)\oplus \Gamma_{mult}(\mathcal{V})$ which is a bijection. It remains to be shown that such a map is a functor. As usual, it suffices to prove this for $\mathcal{V}$ split. For that, consider multiplicative sections $V,V',V''\in \Gamma_{mult}(\mathcal{V})$ covering $e,e',e''\in\Gamma(E)$, respectively. Assume that $\tau:V\Rightarrow V'$ and $\sigma:V'\Rightarrow V''$ are composable morphisms in $\mathrm{Sec}(\G,\mathcal{V})$ with $\tau=(c_0,e_0), \sigma=(c_1,e_1)\in \Gamma(C)\oplus \Gamma(E)$. For every $g\in\G$, there are natural diagrams

\begin{equation}
\xymatrix{
e_{s(g)}   \ar[d]_{\tau(s(g))}  \ar[r]^{V(g)}  &   e_{t(g)}   \ar[d]^{\tau(t(g))}  \\
e'_{s(g)}   \ar[r]_{V'(g)}  &   e'_{t(g)}
}
\end{equation}

\noindent and 

\begin{equation}
\xymatrix{
e'_{s(g)}   \ar[d]_{\sigma(s(g))}  \ar[r]^{V'(g)}  &   e'_{t(g)}   \ar[d]^{\sigma(t(g))}  \\
e''_{s(g)}   \ar[r]_{V''(g)}  &   e''_{t(g)}.
}
\end{equation}

The composition on the groupoid $\mathcal{V}\rr E$ induces a composition of natural transformation, which yields a morphism $\sigma\circ \tau:V\Rightarrow V''$ whose natural diagram at $g\in\G$ is

\begin{equation}\label{functorial naturality}
\xymatrix{
e_{s(g)}   \ar[d]_{\sigma(s(g))\tau(s(g))}  \ar[r]^{V(g)}  &   e_{t(g)}   \ar[d]^{\sigma(t(g))\tau(t(g))}  \\
e''_{s(g)}   \ar[r]_{V''(g)}  &   e''_{t(g)}
}
\end{equation}

We will show that the arrow in $\Gamma(C)\oplus \Gamma_{mult}(\mathcal{V})\rr \Gamma_{mult}(\mathcal{V})$ corresponding to $\sigma\circ \tau: V\Rightarrow V''$ is the composition 

$$(c_1,e_1)(c_0,e_0)=(c_1+c_0,e_0),$$
in $\Gamma(C)\oplus \Gamma_{mult}(\mathcal{V})\rr \Gamma_{mult}(\mathcal{V})$. Note that $2.$ in Proposition \ref{explicitmorphism} implies that $e_1=\delta(c_0)+e_0$, hence the composition above makes sense.

Now, due to Proposition \ref{explicitmorphism} we only need to compute the composition of the arrows in \eqref{functorial naturality} and prove that 

\begin{align}\label{eq:func1}
\sigma(t(g))\tau(t(g))V(g)=&(c_1+c_0)^r(g)+V(g)\\
\label{eq:func2}
V''(g)\sigma(s(g))\tau(s(g))=& V''(g)+(c_1+c_0)^l(g).
\end{align}

We only show identity \eqref{eq:func1}, since \eqref{eq:func2} follows in a similar manner. Since $\tau:V\Rightarrow V'$ is a morphism, we have shown in Proposition \ref{explicitmorphism} that 

$$\tau(t(g))V(g)=c^r_0(g)+V(g)=\big{(}c_0(t(g))+c(t(g)),g,e(s(g))\big{)},$$
so, identities \eqref{eq:splitvbgmaps} imply that the left hand side of \eqref{eq:func1} reads

\begin{align*}
\sigma(t(g))\tau(t(g))V(g)=& \big{(}c_1(t(g)),t(g),e_1(t(g))\big{)}\big{(}c_0(t(g))+c(t(g)),g,e(s(g))\big{)}\\
=& \big{(}c_1(t(g))+c_0(t(g)),g,e(s(g))\big{)}\\
=& (c_1+c_0)^r(g)+V(g).
\end{align*}

This shows identity \eqref{eq:func1}. A similar computation shows \eqref{eq:func2}.

\end{proof}

\subsection{Examples}

Let us discuss some important examples.

\begin{example}
\label{example: sections of 2-vector space}
If $V_1 \rightrightarrows V_0$ is a 2-vector space, considered as a $\mathcal{VB}$-groupoid over $\left\{\ast\right\} \rightrightarrows \left\{\ast\right\}$, then $\mathrm{Sec}\left(V_1 \rightrightarrows V_0\right)$ can be identified with $V_1 \rightrightarrows V_0$ itself. Indeed, category of functors from the trivial groupoid $\left\{\ast\right\} \rightrightarrows \left\{\ast\right\}$ to $V_1 \rightrightarrows V_0$ is canonically isomorphic to (the  category underlying) $V_1 \rightrightarrows V_0$. In terms of Theorem \ref{2vect}, the complex of multiplicative sections is exactly the 2-term complex 
\[
t: \mathrm{Ker} (s) \to V_0
\] 
associated to $V_1 \rightrightarrows V_0$ (see Section \ref{subsection: 2-vector spaces}), and the isomorphism 
\[
V_1 \rightrightarrows V_0  \; \cong \; \mathrm{Ker}(s) \oplus V_0 \rightrightarrows V_0 
\]
is exactly that arising from the equivalence (\ref{eqn: equivalence of 2Vect and 2Term}).
\end{example}

\begin{example}\label{multiplicativefields}

Let $T\G\rr TM$ be the tangent groupoid of $\G\rr M$ defined in Example \ref{tangentVB}. In this case, the category of multiplicative sections coincides with the category of multiplicative vector fields $\mathrm{Vect}(\G)$ introduced by Hepworth in \cite{hepworth}. It follows from Theorem \ref{2vect}, that $\mathrm{Vect}(\G)$ has the structure of a 2-vector space
\[
\mathrm{Vect}(\mathcal G) \; \cong\; \Gamma(A) \oplus \mathfrak X_{mult}(\mathcal G) \rightrightarrows \mathfrak X_{mult}(\mathcal G)
\]
corresponding to the 2-term complex of vector spaces 
\begin{align*}
\Gamma(A) & \to \mathfrak{X}_{mult}(\G) \\
a & \mapsto a^r-a^l
\end{align*}
where $A$ denotes the Lie algebroid of $\G$.

\end{example}

\begin{example}\label{multiplicativeforms}

Let $T^*\G\rr A^*$ the cotangent groupoid of $\G\rr M$ defined in Example \ref{cotangentVB}. In this case, the core bundle is $T^*M\rmap M$ and any 1-form $\alpha\in\Omega^1(M)$ induces right and left invariant 1-forms on $\G$ given by
\[
\alpha^r=t^*\omega; \quad \alpha^l=s^*\omega.
\]
\noindent See Section 11.3 in \cite{mackbook} for more details. Hence, the category of multiplicative 1-forms on $\G$ has a 2-vector space structure
\[
\mathrm{Sec}(\mathcal G,T^\ast \mathcal G)
\; \cong \;  \Omega^1 (\mathcal G) \oplus \Gamma^1_{mult}(\mathcal G) \rightrightarrows \Omega^1_{mult}(\mathcal G)
\]
corresponding to the 2-term complex of vector spaces 
\begin{align*}
\Omega^1(M) & \to \Omega^1_{mult}(\G) \\ \alpha& \mapsto t^*\alpha - s^*\alpha
\end{align*}

\end{example}

\begin{example}
\label{example: Sec of repns}
If $\mathcal G \ltimes E \rightrightarrows E$ is the $\mathcal{VB}$-groupoid associated to a representation of a Lie groupoid $\mathcal G \rightrightarrows M$ on a vector bundle $E \to M$ (see Example \ref{transformationVB}) then the core of $\mathcal G \ltimes E$ is $0$, and $\Gamma_{mult}(\mathcal G \ltimes E) \cong \Gamma(E)^\mathcal G$ (Example \ref{example: sections of transformationVB}). Therefore the complex of multiplicative sections is
\begin{align*}
0 & \to \Gamma\left(E\right)^\mathcal G
\end{align*}
and 
\[
\mathrm{Sec}(\mathcal G, \mathcal G \ltimes E) \; \cong \;
\Gamma\left(E\right)^\mathcal G \rightrightarrows \Gamma\left(E\right)^\mathcal G
\]
\end{example}



\subsection{Regular $\vb$-groupoids}

In this section we describe the category of multiplicative sections of the so-called regular $\vb$-groupoids \cite{GM10}. A $\vb$-groupoid $\mathcal{V}\rr E$ is called \textbf{regular} if the core anchor map $\partial:C\rmap E$ has constant rank. 
There are two extreme classes of regular $\vb$-groupoids; a regular $\vb$-groupoid is of \textbf{type 0} if the core anchor map is zero, and we say that it is of \textbf{type 1} if the core anchor map $\partial:C\rmap E$ is an isomorphism of vector bundles. It turns out that any regular $\vb$-groupoid is isomorphic to the direct sum of a $\vb$-groupoid of type 0 with one of type 1. For more details, see \cite{GM10}.

\subsubsection*{Regular $\vb$-groupoids of type 1}

Let $\G\rr M$ be a Lie groupoid and $q_E:E\rmap M$ a vector bundle. The pullback groupoid $q^!_{E}\G\rr E$ has a natural structure of vector bundle over $\G$ given by the direct sum $q^!_{E}\G=s^*E\oplus t^*E$. It is easy to see that $q^!_E\G\rr E$ is a $\vb$-groupoid over $\G\rr M$, with structure maps

\begin{align*}
\tilde{s}(e_1,g,e_2)=&e_2\\
\tilde{t}(e_1,g,e_2)=&e_1\\
(e_1,g_g,e_2)(e_2,g_2,e_3)=&(e_1,g_1g_2,e_3)\\
(e_1,g,e_2)^{-1}=&(e_2,g^{-1},e_1).
\end{align*}

\noindent In this case, the core anchor is $\mathrm{id}:E\rmap E$, hence $q^!_{E}\G$ is a regular $\vb$-groupoid of type one. As shown in \cite{GM10}, any regular $\vb$-groupoid of type 1 is isomorphic to a pullback $q^!_{E}\G\rr E$. In order to describe the complex of multiplicative sections of a regular groupoid of type 1, it suffices to do it for the pullback groupoid $q^!_{E}\G\rr E$.

It is easy to check that any multiplicative section of $q^!_{E}\G$ has the form $\G\rmap q^!_{E}\G;g\mapsto (e_{t(g)},g,e_{s(g)})$, for some section $e\in\Gamma(E)$ and hence the space of multiplicative sections of $q^!_{E}\G$ identifies with $\Gamma(E)$. Also, if $e\in\Gamma(E)$ is a section of the \emph{core bundle} $E$, then the induced right and left invariant sections $e^r,e^l:\G\rmap q^!_{E}\G$ are

$$e^r(g)=(e_{t(g)},g,0); \quad e^l(g)=(0,g,-e_{s(g)}).$$

We conclude that the complex of multiplicative sections of $q^!_{E}\G\rr E$ is the exact complex

\begin{align}\label{eq:complexregulartype1}
\delta:\Gamma(E)&\rmap \Gamma(E)\cong \Gamma_{mult}(q^!_{E}\G)\nonumber \\ &e\mapsto e.
\end{align}

\subsubsection*{Regular $\vb$-groupoids of type 0}

Let $\mathcal{V}\rr E$ be a regular $\vb$-groupoid of type zero over $\G\rr M$. This is equivalent to have a pair of representations $\Delta^C,\Delta^E$ on $C$ and $E$ respectively, together with a normalized element $\Omega\in \Gamma(\G^{(2)},\mathrm{Hom}(s^*E,t^*C))$. Indeed, $\mathcal{V}\rr E$ is isomorphic to the $\vb$-groupoid $t^*C\oplus s^*E\rr E$ with structure maps

\begin{align}\label{eq:splitmapsregularzero}
\tilde{s}(c,g,e)=&e\nonumber \\
\tilde{t}(c,g,e)=&\Delta^E_g(e)\nonumber \\
(c,g,\Delta^E_h(e))(c',h,e)=&(c+\Delta^C_g(c')-\Omega_{g,h}(e),gh,e).
\end{align}

We will use \eqref{eq:splitmapsregularzero} to describe multiplicative sections on $t^*C\oplus s^*E\rr E$. For that, consider $V:\G\rmap t^*C\oplus s^*E$ a multiplicative section covering $e\in\Gamma(E)$. Since $\tilde{s}\circ V=e\circ s$ and $\tilde{t}\circ V=e\circ t$, we conclude that $e\in\Gamma(E)^{\G}$ is an invariant section, i.e. $e_{t(g)}=\Delta^E_ge_{s(g)}$ for every $g\in\G$. Now, for any composable pair $(g,h)\in \G^{(2)}$ the elements $V(g),V(h)\in t^*C\oplus s^*E$ are composable. Let us write $V(g)=(c_{t(g)},g,e_{t(g)}), V(h)=(c'_{t(h)},h,e_{s(h)})$ and $V(gh)=(\tilde{c}_{t(gh)},gh,e_{t(gh)})$ for some core elements $c_{t(g)}\in C_{t(g)}, c'_{t(h)}\in C_{t(h)}$ and $\tilde{c}_{t(gh)}$. The multiplicativity condition $V(gh)=V(g)V(h)$ reads

\begin{align*}
(\tilde{c}_{t(gh)},gh,e_{t(gh)})=&(c_{t(g)},g,e_{t(g)})(c'_{t(h)},h,e_{s(h)})\\
=& (c_{t(g)}+\Delta^C_gc'_{t(h)}-\Omega_{g,h}e_{s(h)},gh,e_{s(h)}),
\end{align*}

\noindent which implies that $c,c'$ and $\tilde{c}$ are constrained by

\begin{equation}\label{eq:constraint}
\tilde{c}_{t(gh)}-c_{t(g)}-\Delta^C_gc'_{t(h)}+\Omega_{g,h}e_{s(h)}=0.
\end{equation}

As a consequence, we conclude that the space of multiplicative sections on $t^*C\oplus s^*E$ is given by the subspace of $\Gamma(t^*C)\oplus \Gamma(E)^{\G}$ defined by the constraint \eqref{eq:constraint}. 
\begin{remark}
It is worthwhile to observe that if $C=0$ the constraint \eqref{eq:constraint} holds automatically and hence $\Gamma_{mult}(s^*E)\cong\Gamma(E)^{\G}$, as seen in Example \ref{example: Sec of repns}.
\end{remark}

Now we are able to describe the complex of multiplicative sections of the type zero $\vb$-groupoid $t^*C\oplus s^*E\rr E$. It follows from \eqref{splitrightinvariant} and \eqref{splitleftinvariant} that the right and left invariant sections associated to $c\in\Gamma(C)$ are respectively given by 

$$c^r(g)=(c_{t(g)},g,0); \quad c^l(g)=(\Delta^C_gc_{s(g)},g,0).$$

Therefore, the complex of multiplicative sections is given by

\begin{align}\label{eq:complexregulartypezero}
\delta:\Gamma(C)\rmap &\Gamma_{mult}(t^*C\oplus s^*E)\\ &c\mapsto c^r-c^l,
\end{align}
where $(c^r-c^l)(g)=(c_{t(g)}-\Delta^C_gc_{s(g)},g,0)$ for every $g\in\G$.

\subsubsection*{General regular $\vb$-groupoids}

Consider now a regular $\vb$-groupoid $\mathcal{V}$ with core anchor $\partial:C\rmap E$. If $K=\ker(\partial), F=\mathrm{im}(\partial)$ and $\nu=E/F=\mathrm{coker}(\partial)$, one can construct two $\vb$-groupoids $\mathcal{V}_0$ and $\mathcal{V}_1$ with $\mathcal{V}\cong \mathcal{V}_0\oplus \mathcal{V}_1$, namely

\begin{equation}\label{eq:decompositionregularVB}
\mathcal{V}_0:=t^*K\oplus s^*\nu\rr \nu; \quad \mathcal{V}_1:=t^*F\oplus s^*F\rr F.
\end{equation}

Note that $\mathcal{V}_0$ has core anchor $\partial:K\rmap \nu$, hence $\mathcal{V}_0$ is regular of type zero. Similarly, the core anchor of $\mathcal{V}_1$ is $id:F\rmap F$, hence $\mathcal{V}_1$ is regular of type one. For more details the reader can consult \cite{GM10}. 

As a result, we obtain that the complex of multiplicative sections of any regular $\vb$-groupoid with core $C\cong K\oplus F$ and side $E\cong F\oplus \nu$, is isomorphic to the complex

\begin{align}\label{eq:complexregular}
\delta:\Gamma(K)\oplus \Gamma(F)&\rmap \Gamma_{mult}(t^*K\oplus s^*\nu)\oplus \Gamma(F)\\
\label{eq:complexregular1}
&k\oplus v\rmap (k^r-k^l)\oplus v, 
\end{align}

\noindent with $k^r-k^l$ defined as in \eqref{eq:complexregulartypezero}.

\begin{corollary}
The complex $C^\bullet_{mult}(\mathcal V)$ of multiplicative sections of a regular $\mathcal{VB}$-groupoid $\mathcal V$ is quasi-isomorphic to the complex $C^\bullet_{mult}(\mathcal V_0)$ of multiplicative sections of its type zero component.
\end{corollary}

\begin{proof}
It is clear from (\ref{eq:complexregular1}) that the inclusion of the subcomplex $\delta:\Gamma(K) \to \Gamma_{mult}(t^*K\oplus s^*\nu)$ into (\ref{eq:complexregular}) is a quasi-isomorphism.
\end{proof}

\begin{remark} We will see later (Corollary \ref{regularqisotypezero}) that this result can also be seen as a consequence of Morita invariance.
\end{remark}

\begin{example}\label{ex:complexregularTG}
A Lie groupoid $\G\rr M$ is called \textbf{regular} if $t:s^{-1}(x)\rmap M$ has constant rank for every $x\in M$. Infinitesimally, the anchor map $\rho:A\rmap TM$ has constant rank. In particular, the tangent $\vb$-groupoid $T\G\rr TM$ is regular. By Example \ref{multiplicativefields}, we know that the complex of multiplicative sections of $T\G$ is just the complex $\Gamma(A)\rmap \mathfrak{X}_{mult}(\G)$ of multiplicative vector fields on $\G$. Hence, the choice of a splitting of $T\G\rr TM$ induces an isomorphism of complexes 

$$\Big{(}\Gamma(A)\rmap \mathfrak{X}_{mult}(\G)\Big{)}\cong \Big{(}\Gamma(K)\oplus \Gamma(F)\rmap \Gamma_{mult}(t^*K\oplus s^*\nu)\oplus \Gamma(F)\Big{)},$$

\noindent where $K=\ker(\rho), F=\mathrm{im}(\rho)$ and $\nu=TM/F$. This example will be studied in detail in Section \ref{sec:vectorfields}.

\end{example}



\section{Cohomology}\label{sec:cohomology}

In this section we study the cohomology of the complex of multiplicative sections as well as its connection with the theory of representations up to homotopy.

\subsection{Cohomology of $\mathcal{VB}$-groupoids}


Let $\G\rr M$ be a Lie groupoid. The \textbf{nerve} of $\G$ is the simplicial manifold $\G_{\bullet}=(\G_{(p)})_{p\in \N}$, where $\G_{(p)}:=\{(g_1,...,g_p)\in\G\times...\times\G;s(g_i)=t(g_{i+1}), i=1,...,p-1\}$ is the manifold of composable $p$-arrows and face maps $\partial_i:\G_{(p)}\to \G_{(p-1)} , i=0,...,p$, defined by

$$\partial_i(g_1,...,g_p)=\begin{cases}(g_2,...,g_p); \quad i=0\\ (g_1,...,g_ig_{i+1},...,g_p); \quad i=1,...,p-1\\ (g_1,...,g_{p-1}); \quad i=p\end{cases}.$$

As for any simplicial manifold, one can define a complex $(C^{\bullet}(\G),d)$ where $d:C^{\infty}(\G_{(p)})\rmap C^{\infty}(\G_{(p+1)})$ is given by

\begin{equation}\label{gpddifferential}
\delta f(g_1,...,g_{p+1})=\sum_{i=0}^{p}(-1)^i\partial^*_i.
\end{equation}

\noindent The corresponding cohomology is called the \textbf{groupoid cohomology} of $\G$ and it is denoted by $H^{\bullet}(\G)$.

Consider now a $\vb$-groupoid $\mathcal{V}\rr E$ over $\G\rr M$ with core $C$. It was proved in \cite{cabreradrummond} that the nerve $\mathcal{V}_{\bullet}$ is a simplicial vector bundle over the nerve $\G_{\bullet}$. The linear complex $(C^{\bullet}_{lin}(\mathcal{V}),d)$ of a $\vb$-groupoid was introduced in \cite{cabreradrummond}. A \textbf{linear $p$-cochain} is a smooth function $\mathcal{V}_{(p)}\rmap \R$ which is fiberwise linear with respect to the vector bundle structure $\mathcal{V}_{(p)}\rmap \G_{(p)}$.  The differential $d$ is just the restriction of the groupoid differential \eqref{gpddifferential}. It was also proved in \cite{cabreradrummond} that $(C^{\bullet}_{lin}(\mathcal{V}),d)$ is a subcomplex of $(C^{\bullet}(\mathcal{V}),d)$, whose cohomology is denoted by $H^{\bullet}_{lin}(\mathcal{V})$.

\subsection{Multiplicative sections and linear cochains}

Let $E \to M$ be a vector bundle and $C^\infty_{lin}\left(E^\ast\right)$ denote the vector space of fibre-wise linear smooth functions on $E^\ast$. There is a linear isomorphism
\begin{align*}
l : \Gamma(E) & \to C^\infty_{lin}\left(E^\ast\right) \\
e & \mapsto \left( \alpha \in E^\ast_m \mapsto \alpha \left( e(m)\right) \right) 
\end{align*}
from $\Gamma(E)$ to the vector space of fibre-wise linear smooth functions on $E^\ast$. (We will denote this map by $l$ for any vector bundle.) These maps satisfy several naturality properties, see Proposition \ref{prop: naturality of l}.

Recall that a $\vb$-groupoid $\mathcal{V}\rr E$ over $\G\rr M$ with core $C$ has a dual $\vb$-groupoid $\mathcal{V}^*\rr C^*$ over $\G\rr M$ with core $E^*$ (c.f. Section 11.2 in \cite{mackbook}).  In particular, any section $V:\G\rmap \mathcal{V}$ induces a linear 1-cochain $l_V$ on the dual $\vb$-groupoid $\mathcal{V}^*\rr C^*$. 

\begin{proposition}\label{mult1cocycles}

The canonical isomorphism $l: \Gamma(\mathcal{V}) \to C^{\infty}_{lin}(\mathcal{V}^*)$ restricts to an isomorphism between multiplicative sections of $\mathcal{V}$ and linear 1-cocycles on $\mathcal{V}^*$

$$l : \Gamma_{mult}(\mathcal{V})\xrightarrow{\cong} Z^1_{lin}(\mathcal{V}^*)$$

\end{proposition}

\begin{proof}

Let $l_V:\mathcal{V}^*\rmap \R$ be the linear 1-cochain associated to a section $V:\G\rmap \mathcal{V}$. Then $l_V$ is a cocycle if and only if $l_V(\xi_g\eta_h)=l_V(\xi_g)+l_V(\eta_h)$ for any pair of composable arrows $g,h\in\G$. One observes that $l_V(\xi_g\eta_h)=\xi_g\eta_h(V(gh))$ and $l_V(\xi_g)+l_V(\eta_h)=\xi_g(V(g))+\eta_h(V(h))$. Hence, $l_V$ is a cocycle if and only if $\xi_g\eta_h(V(gh))=\xi_g(V(g))+\eta_h(V(h))$, which is equivalent to saying that $V$ is multiplicative.
\end{proof}

\begin{proposition}\label{1truncationisomorphism}

The linear isomorphisms $l : \Gamma(C) \to C^\infty_{lin}\left(C^\ast\right)$ and $l : \Gamma_{mult}\left(\mathcal V\right) \to Z^1_{lin}\left(\mathcal V^\ast\right)$ define an isomorphism of 2-term complexes 
\[
C^\bullet_{mult} \left(\mathcal V\right) \to
C^\bullet_{lin}\left(\mathcal V^\ast\right)^{\leq 1}
\]
where $C^\bullet_{lin}\left(\mathcal V^\ast\right)^{\leq 1}$ is the 1-truncation $d:C^{\infty}_{lin}(C^*)\rmap Z^1_{lin}(\mathcal{V}^*)$ of $C^\bullet_{lin}\left(\mathcal V^\ast\right)$.
\end{proposition}

\begin{proof}
It follows from Proposition \ref{mult1cocycles} that the statement is equivalent to the commutativity of the diagram
\[
\xymatrix{
C^\infty_{lin}\left(C^\ast\right)  \ar[r]^d  &  Z^1_{lin}\left(\mathcal V^\ast\right)  \\
\Gamma(C)  \ar[u]^l  \ar[r]_\delta  &  \Gamma_{mult}\left(\mathcal V\right) \ar[u]_l
}
\]
Denote the structure maps of $\mathcal V^\ast \rightrightarrows C^\ast$ by $\hat s, \hat t, \hat m, \hat i, \hat 1$, and the zero section by $\hat 0: \mathcal G \to \mathcal V^\ast$. Recall from \cite{mackbook} that if $\alpha \in \mathcal V^\ast_g$ then $\hat t(\alpha) \in C^\ast_{t(g)}$, $\hat s(\alpha) \in C^\ast_{s(g)}$, and $\hat s(\alpha)$ and $\hat t(\alpha)$ are given by
\begin{align*}
\hat s (\alpha) (u) & = \alpha \left(  \hat 0_g \cdot \hat i u\right) \\
\hat t (\alpha) (u) &  = \alpha \left( u \cdot \hat 0_g\right)
\end{align*}
It follows that if $c \in \Gamma(C)$ then
\begin{align*}
\left(d l_c  \right) (\alpha)& = \left(\hat t^\ast \left(l_c\right)- \hat s^\ast \left(l_c\right) \right)(\alpha) \\
& = \left( l_c \hat t-  l_c \hat s\right)(\alpha) \\
& = l_c \left( \hat t (\alpha) \right) - l_c \left(\hat s(\alpha) \right) \\
& = \hat t (\alpha)  \left(c(t(g)\right) 
- \hat s (\alpha)  \left(c(s(g)\right) \\
& = \alpha \left(   c(t(g)) \cdot  \hat 0_g \right) 
- \alpha \left( \hat 0_g \cdot \hat i c(s(g)) \right) \\
& = \alpha \left(   c(t(g)) \cdot  \hat 0_g -  \hat 0_g \cdot \hat i c(s(g)) \right) \\
& = \alpha \left(  \left(c^r - c^l\right)(g)  \right) \\
& = \alpha \left(  \left(\delta c\right)(g)  \right) \\
& = l_{\delta c} \left(\alpha\right)
\end{align*}
Therefore $dl = l \delta$ as required.
\end{proof}

\subsection{Representations up to homotopy}

As shown in \cite{cabreradrummond}, the cohomology of the linear complex $C^{\bullet}_{lin}(\mathcal{V}^*)$ is isomorphic to $H^{\bullet}(\G,\mathcal{E}_{\mathcal{V}})$, the cohomology of $\G$ with coefficients in the representation up to homotopy $\mathcal{E}_{\mathcal{V}}$ corresponding to the $\vb$-groupoid $\mathcal{V}$. In particular, we obtain the following result which describes the cohomology of the category of multiplicative sections $\mathrm{Sec}(\G,\mathcal{V})$. 

\begin{proposition}\label{cohomology}

The cohomology of the 2-term complex of multiplicative sections $\delta:\Gamma(C)\rmap \Gamma_{mult}(\mathcal{V})$ is given by $\Gamma(C)^\G:=\{c\in\Gamma(C);c^r=c^l\}$ in degree 0 and $H^1(\G,\mathcal{E}_{\mathcal{V}})$ in degree 1.

\end{proposition}

A result of Arias Abad and Crainic \cite{AriasCrainic2}, establishes that if $\G$ is proper Lie groupoid and $\mathcal{E}$ is a representation up to homotopy of $\G$, then the cohomology of $\G$ with coefficients in $\mathcal{E}$ is concentrated in degree 0 and degree 1. Combining this result with Proposition \ref{cohomology}, we obtain the following.

\begin{corollary}\label{cor:cohomologyruth}

Let $\mathcal{V}\rr E$ be a proper $\vb$-groupoid. Then the 2-term complex of multiplicative sections $\delta:\Gamma(C)\rmap \Gamma_{mult}(\mathcal{V})$ computes the cohomology of $\G$ with coefficients in the representation up to homotopy $\mathcal{E}_{\mathcal{V}}$ corresponding to $\mathcal{V}$.

\end{corollary}

\begin{proof}
One easily observes that if $\mathcal{V}$ is proper, then $\G$ is proper as well. Hence, the cohomology $H^{\bullet}(\G,\mathcal{E}_{\mathcal{V}})$ is concentrated in degrees 0 and 1, which can be computed by using the complex of multiplicative sections of $\mathcal{V}$. 
\end{proof}

\begin{remark}
Notice that a $\vb$-groupoid over a proper groupoid $\G\rr M$ is not necessarily a proper $\vb$-groupoid. For instance, any 2-vector space is a $\vb$-groupoid over the proper groupoid $*\rr *$. However, 2-vector spaces are not necessarily proper.
\end{remark}









\section{Morita invariance}\label{sec:moritainvariance}

In this section we study Morita invariance of the category of multiplicative sections $\mathrm{Sec}(\G,\mathcal{V})$ of a $\vb$-groupoid $\mathcal{V}$ as in \eqref{VBgroupoid}. We first explain the notion of Morita equivalence of $\vb$-groupoids, introduced in \cite{dHO}.

\subsection{$\mathcal{VB}$-Morita maps}

 Recall that a groupoid morphism $\varphi:\G_1 \to \G_2$ covering $f:M_1\rmap M_2$ is called a \textbf{weak equivalence} if it satisfies the following conditions:

\begin{enumerate}

\item the map $\G_{2}\times _{M_2}^{s_2,f}M_1\rmap M_2; (g,x)\mapsto t(g)$ is a surjective submersion
\item the square

$$\xymatrix{\ar@{}[dr]|{}
\G_1 \ar[r]^\varphi \ar[d]_{(s\times  t)}
& \G_2 \ar[d]^{(s,t)}\\
M_1\times  M_1 \ar[r]^{f\times f} & M_2\times M_2,}$$

\noindent is a good fiber product of manifolds \cite{dH}.

\end{enumerate}

Following \cite{dHO}, we specialize the notion of Morita map to the framework of $\vb$-groupoids.

\begin{definition}
\label{def: VB Morita maps}

Let $\mathcal{V}_1$ and $\mathcal{V}_2$ be $\vb$-groupoids over $\G_1$ and $\G_2$, respectively. A morphism of $\mathcal{VB}$-groupoids $\left(\Phi,\phi\right):\mathcal V_1 \to \mathcal V_2$ is called a \textbf{$\mathcal{VB}$-Morita map} if $\Phi$ is a weak equivalence of Lie groupoids.
\end{definition}

For more details about this notion, the reader can see \cite{dHO}. Following the case of Lie groupoids, the notion of $\mathcal{VB}$-Morita maps leads to a notion of Morita equivalence for $\mathcal{VB}$-groupoids:

\begin{definition}
We say that $\mathcal{VB}$-groupoids $\mathcal V$ and $\mathcal V'$ are \emph{Morita equivalent} if there exists a third $\mathcal{VB}$-groupoid $\mathcal W$ and $\mathcal{VB}$-Morita morphisms $\Phi:\mathcal W \to \mathcal V$ and $\Phi':\mathcal W \to \mathcal V'$ as in the following diagram:
\[
\xymatrix{
& \mathcal W \ar[dl]_{\Phi'}  \ar[dr]^{\Phi'} & \\
\mathcal V  & & \mathcal V'
}
\]
\end{definition}

\begin{proposition} 
\label{prop: VB Morita is an equivalence relation}
Morita equivalence for $\mathcal{VB}$-groupoids is an equivalence relation.
\end{proposition}

\begin{proof}
Using Proposition \ref{prop: existence of VB weak fibre products} the proof is identical to the analogous statement for Lie groupoids (\cite{MoerdijkMrcun}).
\end{proof}

\begin{theorem}\cite{dHO}\label{VBmoritacohomology} Let $\Phi:\mathcal{V}_1\rmap \mathcal{V}_2$ be a $\vb$-Morita map. The chain map $\Phi^*:C^{\bullet}_{lin}(\mathcal{V}_2)\rmap C^{\bullet}_{lin}(\mathcal{V}_1)$ is a quasi-isomorphism, i.e. $\Phi^\ast$ induces an isomorphism in cohomology $\Phi^*:H^{\bullet}_{lin}(\mathcal{V}_2)\rmap H^{\bullet}_{lin}(\mathcal{V}_1).$

\end{theorem}

\begin{corollary}
\label{cor: Morita invariance of VB cohomology}
If $\mathcal V$ and $\mathcal V'$ are Morita equivalent $\mathcal{VB}$-groupoids then $H^{\bullet}_{lin}(\mathcal{V}) \cong H^{\bullet}_{lin}(\mathcal{V'})$.
\end{corollary}

\begin{proof}
By Theorem \ref{VBmoritacohomology} a  pair of $\mathcal{VB}$ -Morita morphisms $\Phi:\mathcal W \to \mathcal V$ and $\Phi':\mathcal W \to \mathcal V'$ induces a pair of quasi-isomorphisms
\[
\xymatrix{
C^\bullet_{lin}(\mathcal V) \ar[r]^-{\Phi^\ast}  &
C^\bullet_{lin}(\mathcal W)  & 
\ar[l]_-{{\Phi'}^\ast} C^\bullet_{lin}(\mathcal V')
}
\]
and therefore an isomorphism 
\[
\xymatrixcolsep{4pc}
\xymatrix{
H^\bullet_{lin}(\mathcal V) \ar[r]^-{\left({\Phi'}^\ast\right)^{-1} \circ \Phi^\ast} & H^\bullet_{lin}(\mathcal V')
}
\]
\end{proof}

\subsection{Notation}
\label{subsection: VB notation}

For the remainder of Section \ref{sec:moritainvariance} let $\mathcal V \rightrightarrows E$ and $\mathcal W \rightrightarrows F$ be $\mathcal{VB}$-groupoids with bases $\mathcal G \rightrightarrows M$ and $\mathcal H \rightrightarrows N$ respectively. Denote the structure maps of $\mathcal V,\mathcal G, \mathcal W, \mathcal H$ by $\tilde s,...,s,...,\tilde s',...,s',...$ respectively. Denote the cores of $\mathcal V$ and $\mathcal W$ by $C$ and $D$ respectively, and the differentials of the complexes $C_{mult}^\bullet(\mathcal V)$ and $C_{mult}^\bullet(\mathcal W)$ by $\delta_\mathcal V$ and $\delta_\mathcal W$, or just by $\delta$ if no confusion will arise. We also fix a morphism of $\mathcal{VB}$-groupoids $\Phi: \mathcal V \to \mathcal W$ covering a groupoid morphism $\phi:\mathcal G \to \mathcal H$:
\[
\xymatrix{
\mathcal V  \ar[d]  \ar[r]^{\Phi} &  \mathcal W \ar[d] \\
\mathcal G  \ar[r]_{\phi}  &  \mathcal H
}
\]
Denote by $\Phi_C : C \to D$ the map determined by $\Phi$. This is a vector bundle morphism covering $\phi_0 : M \to N$:
\[
\xymatrix{
C  \ar[d]  \ar[r]^{\Phi_C} &  D \ar[d] \\
M  \ar[r]_{\phi_0}  &  N
}
\]
Our aim is now to prove the analogue of Corollary \ref{cor: Morita invariance of VB cohomology} for complexes of multiplicative sections. This will follow from Proposition \ref{1truncationisomorphism}, Theorem \ref{VBmoritacohomology}, and the use of the pullback construction for $\mathcal{VB}$-groupoids. 

\subsection{Pullbacks of multiplicative sections}

As explained in \cite{BursztynCabreradelHoyo}, the pullback vector bundle $\phi^\ast \mathcal W$ is a $\mathcal{VB}$-groupoid over $\mathcal G$, with groupoid structure given by the diffeomorphism 
\[
\phi^\ast \mathcal W \cong \mathcal G \times_\mathcal H \mathcal W
\]
where $\mathcal G\times_\mathcal H \mathcal W$ is the (strict) fibre product (in the category of Lie groupoids) of $\phi: \mathcal G \to \mathcal H$ and the projection $\mathcal W \to \mathcal H$. The projection $\phi^\ast \mathcal W \to \mathcal W$ is a morphism of $\mathcal{VB}$-groupoids, denoted $\phi^!$, and there is a unique morphism $\bar{\Phi}: \mathcal V \to \mathcal \phi^\ast \mathcal W$ such that $\Phi = \phi^! \circ \bar \Phi$:
\begin{equation}
\label{VBmappullback1}
\xymatrix{
\mathcal V  \ar[r]^{\bar \Phi}   \ar@/^{2pc}/[rr]^{\Phi} \ar[d]  &
\phi^\ast \mathcal W \ar[d]  \ar[r]^{\phi^!} &
\mathcal W \ar[d] \\
\mathcal G \ar@{=}[r] & \mathcal G \ar[r]_{\phi} & \mathcal H
}
\end{equation}
It follows from the form of the source morphism of $\phi^\ast \mathcal W$ 
\begin{align*}
\tilde s'' : \mathcal G \times_\mathcal H \mathcal W &  \to M \times_N F \\
\left(g,w\right) & \mapsto \left(s(g),\tilde s'(w) \right)
\end{align*}
that the core of $\phi^\ast \mathcal W$ is isomorphic (as a vector bundle over $\mathcal G$) to $\phi_0 ^\ast D$. In particular, there is a canonical linear map
\begin{align}
\label{eqn: pullback of core sections}
\phi_0^\ast : \Gamma(D) & \to \Gamma\left(\phi_0^\ast D\right) \cong \mathrm{Core} \left(\phi^\ast \mathcal W\right)
\end{align}
As $\phi^\ast \mathcal W$ is a strict fibre product of Lie groupoids, any multiplicative section $V$ of $\mathcal W$ pulls back to a multiplicative section $\phi^\ast V$ of $\phi^\ast \mathcal W$, which induces a linear map 
\begin{align}
\label{eqn: pullback of mult sections}
\phi^\ast : \Gamma_{mult}(\mathcal W) & \to \Gamma_{mult}\left(\mathcal \phi^\ast \mathcal W\right)
\end{align}
\begin{proposition}
\label{prop: pullback chain map}
The linear maps (\ref{eqn: pullback of core sections}) and (\ref{eqn: pullback of mult sections}) are the components of a chain map
\begin{equation}
\label{eqn: pullback chain map}
\phi^\ast : C^\bullet_{mult}\left(\mathcal W\right) \to
C^\bullet_{mult}\left(\phi^\ast\mathcal W\right)
\end{equation}
\end{proposition}

\begin{proof}
If $d \in \Gamma(D)$ then
\begin{align*}
\left(\phi_0^\ast d\right)^r & = \left((\phi_0^\ast d)(t(g))\right) \tilde 0_g \\
& = \left( 1_{t(g)}, d(\phi_0 t(g)) \right)  \left( g,  \tilde 0_{\phi(g)}\right) \\
& = \left( 1_{t(g)}, d(t\phi(g)) \right)  \left( g,  \tilde 0_{\phi(g)}\right) \\
& = \left( g , d(t\phi(g)) \tilde 0_{\phi(g)}\right) \\
& = \left( g, d^r (\phi(g)) \right) \\
& = \left(\phi ^\ast d^r\right)(g)
\end{align*}
and therefore
\begin{align*}
\left( \phi_0^\ast d\right)^r & = \phi^\ast \left(d^r\right)
\end{align*}
A similar computation shows that
\begin{align*}
\left( \phi_0^\ast d\right)^l & = \phi^\ast \left(d^l\right)
\end{align*}
and therefore that
\begin{equation*}
\delta_{\phi^\ast \mathcal W} \left(\phi_0^\ast d\right) 
= \left(\phi_0^\ast d\right)^r - \left(\phi_0^\ast d\right)^l
= \phi^\ast \left(d^r\right)  - \phi^\ast \left(d^l\right)
= \phi^\ast \left(\delta_\mathcal W d\right)
\end{equation*}
which shows that
\[
\delta_{\phi^\ast \mathcal W} \phi_0 ^\ast 
=
\phi^\ast \delta_\mathcal W
\]
as required.
\end{proof}

As $\bar \Phi$ is a base-preserving morphism of $\mathcal{VB}$-groupoids it induces a chain map
\begin{equation}
\label{eqn: pullback chain map 2}
\bar \Phi : C_{mult}^\bullet \left(\mathcal V\right)  \to C_{mult}^\bullet\left(\phi^\ast \mathcal W\right)
\end{equation}
We will show in the next section that if $\Phi$ is a $\mathcal{VB}$-Morita map then the chain maps (\ref{eqn: pullback chain map}) and (\ref{eqn: pullback chain map 2}) are quasi-isomorphisms.

\subsection{Morita invariance of $C^\bullet_{mult}(\mathcal V)$}

If we pullback the dual $\mathcal{VB}$-groupoid $\mathcal W^\ast$ along $\phi:\mathcal G \to \mathcal H$, then after identifying $\phi^\ast \mathcal W$ with $\left(\phi^\ast \mathcal W\right)^\ast$, the dual $\bar \Phi^\ast$ of $\bar \Phi$ is a morphism $\bar \Phi^\ast: \phi^\ast \mathcal W^\ast \to \mathcal V^\ast$, and we have a pair of morphisms:
\begin{equation}
\label{VBmappullback2}
\xymatrix{
\mathcal V^\ast \ar[d] &\ar[l]_{\bar{\Phi}^\ast} \phi^\ast \mathcal W^\ast  \ar[d] \ar[r]^-{\phi^!} & \mathcal W^\ast\ar[d] \\
\mathcal G \ar@{=}[r] &  \mathcal G \ar[r]_{\phi} & \mathcal H
}
\end{equation}

\begin{proposition}
\label{prop: pair of VB Morita maps}
If $\Phi:\mathcal V \rmap \mathcal W$ is a $\mathcal{VB}$-Morita map then the morphisms $\bar \Phi^\ast$ and $\phi^!$ in (\ref{VBmappullback2}) are $\mathcal{VB}$-Morita maps.
\end{proposition}

\begin{proof}
It is shown in \cite{dHO} that $\phi^! : \phi^\ast \mathcal W ^\ast \to \mathcal W^\ast$ is a Morita map whenever $\phi$ is a weak equivalence, and so in particular whenever $\Phi$ is a Morita map. It then follows from the 2 out of 3 property of weak equivalences that $\bar \Phi$ is a Morita map whenever $\Phi$ is. Finally, it is shown in \cite{dHO} that a `base preserving' morphism of $\mathcal{VB}$-groupoids is a Morita map if and only if its dual is, which implies that $\bar \Phi^\ast$ is a Morita map whenever $\Phi$ is.
\end{proof}

\begin{proposition}
\label{prop: VB qisoms}
If $\Phi$ is a $\mathcal{VB}$-Morita map then the chain maps
\[
\xymatrix{
C_{mult}^\bullet(\mathcal V)  \ar[r]^-{\bar\Phi}  &   C_{mult}^\bullet (f^\ast \mathcal W)     &
\ar[l]_-{\phi^\ast} C_{mult}^\bullet(\mathcal W) 
}
\]
((\ref{eqn: pullback of core sections}) and (\ref{eqn: pullback of mult sections}) above) are quasi-isomorphisms.
\end{proposition}

\begin{proof}
Consider the following diagram of chain maps
\begin{equation}
\label{diag: pullbacks and sections}
\xymatrix{
C_{lin}^\bullet(\mathcal V^\ast)^{\leq 1}  \ar[r]^{\left(\bar \Phi^\ast\right)^\ast} &  C_{lin}^\bullet(f^\ast \mathcal W^\ast)^{\leq 1}  &  \ar[l]_-{\left(\phi^!\right)^\ast} C_{lin}^\bullet(\mathcal W^\ast)^{\leq 1} \\
C_{mult}^\bullet(\mathcal V)  \ar[u]^l \ar[r]_-{\bar\Phi}  &   C_{mult}^\bullet (f^\ast \mathcal W)  \ar[u]_l   &
\ar[l]^-{\phi^\ast} C_{mult}^\bullet(\mathcal W) \ar[u]_l
}
\end{equation}
in which the vertical arrows are the isomorphisms of Proposition \ref{1truncationisomorphism}, and the morphisms in the top row are the 1-truncations of the chain maps induced by the two morphisms of $\mathcal{VB}$-groupoids in (\ref{VBmappullback2}). It follows from Proposition \ref{prop: pair of VB Morita maps}, Theorem  \ref{VBmoritacohomology} and the fact that the 1-truncation functor preserves quasi-isomorphisms that the morphisms in the top row of (\ref{diag: pullbacks and sections}) are both quasi-isomorphisms. Therefore, it is sufficient to show that the diagram (\ref{diag: pullbacks and sections}) commutes, which follows from Proposition \ref{prop: naturality of l}.
\end{proof}

\begin{theorem} \label{moritainvarianceVB}
If $\Phi:\mathcal V \rmap \mathcal W$ is a $\mathcal{VB}$-Morita map covering $\phi:\G \rmap \mathcal H$ then there is an induced isomorphism $H^\bullet (C^\bullet_{mult}(\mathcal V)) \cong H^\bullet (C^\bullet_{mult}(\mathcal W))$.
\end{theorem}

\begin{proof}
The quasi-isomorphisms of Proposition \ref{prop: VB qisoms} induce isomorphisms between homology groups
\[
\xymatrix{
H^\bullet\left(C^\bullet_{mult}(\mathcal V)\right) \ar[r]^-{\bar\Phi}  &
H^\bullet\left(C^\bullet_{mult}(\mathcal \phi^\ast \mathcal W)\right)  & 
\ar[l]_-{\phi^\ast} 
H^\bullet \left(C^\bullet_{mult}(\mathcal W)\right)
}
\]
and therefore an isomorphism 
\[
\xymatrixcolsep{4pc}
\xymatrix{
H^\bullet\left(C^\bullet_{mult}(\mathcal V)\right) \ar[r]^-{\left(\phi^\ast\right)^{-1} \circ \bar \Phi} & H^\bullet \left(C^\bullet_{mult}(\mathcal W)\right)
}
\]
\end{proof}

More generally, we have:

\begin{theorem} 
\label{moritainvarianceVB2}
If $\mathcal V$ and $\mathcal V'$ are Morita equivalent $\mathcal{VB}$-groupoids then $H^\bullet (C^\bullet_{mult}(\mathcal V)) \cong H^\bullet (C^\bullet_{mult}(\mathcal V'))$.
\end{theorem}

\begin{proof}
If $\mathcal W \to \mathcal V$ and $\mathcal W \to \mathcal V'$ are $\mathcal{VB}$-Morita maps then by Theorem \ref{moritainvarianceVB} there are isomorphisms
\begin{equation*}
H^\bullet (C^\bullet_{mult}(\mathcal V)) \cong 
H^\bullet (C^\bullet_{mult}(\mathcal W)) \cong
H^\bullet (C^\bullet_{mult}(\mathcal V'))
\end{equation*}
\end{proof}

\begin{remark}
Note that the isomorphism of Theorem \ref{moritainvarianceVB2} depends on the choice of a $\mathcal{VB}$-groupoid $\mathcal W$ and $\mathcal{VB}$-Morita maps $\mathcal W \to \mathcal V$ and $\mathcal W \to \mathcal V'$.
\end{remark}

\begin{remark}
The proof of Theorem \ref{moritainvarianceVB} actually shows that there is a `zig-zag' of quasi-isomorphisms relating the complexes $C^\bullet_{lin}(\mathcal V)$ and $C^\bullet_{lin}(\mathcal V')$. However, every cochain complex of vector spaces is homotopy equivalent to its cohomology, and so the existence of such a zig-zag is equivalent to the existence of an isomorphism $H^{\bullet}_{lin}(\mathcal{V}) \cong H^{\bullet}_{lin}(\mathcal{V'})$. This is no longer the case in the category of differential graded Lie algebras - in particular in the category of Lie 2-algebras. For this reason the statements of Theorems \ref{theorem: invariance of Lie 2-algebras} and \ref{theorem: L infinity map} (see below) - which are $\mathcal{LA}$-groupoid analogues of Theorems \ref{moritainvarianceVB} and \ref{moritainvarianceVB} above, involve either  the derived category of crossed modules or $L_\infty$ morphisms.
\end{remark}

\subsection{Consequences}

Recall that the complex of multiplicative sections of a regular $\vb$-groupoid $\mathcal{V}\rr E$ over $\G\rr M$ was described in \eqref{eq:complexregular} in terms of the corresponding complexes associated to the type zero and type one components of $\mathcal{V}$. We finish this section with a simple application of Theorem \ref{moritainvarianceVB} which describes the complex of multiplicative sections of a regular $\vb$-groupoid only in terms of its type zero component. For that, consider a regular $\vb$-groupoid $\mathcal{V}\rr E$ over $\G\rr M$ with core anchor $\partial:C\rmap E$. Let $K=\ker(\partial)$, $F=\mathrm{im}(\partial)$ and $\nu=\mathrm{coker}(\partial)=E/F$. Consider the decomposition of $\mathcal{V}\cong \mathcal{V}_0\oplus \mathcal{V}_1$ as in \eqref{eq:decompositionregularVB}.

\begin{corollary}\label{regularqisotypezero}
The complex of multiplicative sections of $\mathcal{V}$ is quasi-isomorphic to the complex of multiplicative sections of its type zero component $\mathcal{V}_0$.

\end{corollary}

\begin{proof}
It suffices to show that the canonical projection $\Phi:\mathcal{V}\rmap \mathcal{V}_0$ is a $\vb$-Morita map, hence the result follows by Theorem \ref{moritainvarianceVB}. In order to prove that $\Phi$ is a $\vb$-Morita map, we will use a result of \cite{dHO} which says that a morphism of $\vb$-groupoids is a $\vb$-Morita map if and only if it is a Morita map between the base groupoids and a quasi-isomorphism fiberwise. In our case, the base map is the identity, hence a Morita map, and clearly $\Phi:\mathcal{V}\rmap \mathcal{V}_0$ induces a quasi-isomorphism between the fibers.

\end{proof}

\section{Projectable sections of $\mathcal{VB}$-groupoids}

In this section we study the sections of a $\mathcal{VB}$-groupoid $\mathcal V$ that are projectable with respect to a morphism $\Phi:\mathcal V \to \mathcal W$. The motivation is to relate the Lie 2-algebras of sections of Morita equivalent $\mathcal{LA}$-gropoids - see Section \ref{section: sections of LA}. (The issue is that if $\mathcal V$ and $\mathcal W$ are $\mathcal{LA}$-groupoids then the $\mathcal{VB}$-groupoid $f^\ast \mathcal W$ in Proposition \ref{prop: VB qisoms} has no natural Lie algebroid structure. Instead, we study projectable sections of $\mathcal V$, and their projection to sections of $\mathcal W$ - an operation that respects the Lie brackets.)

\subsection{Projectable sections of vector bundles}

In this subsection we fix some notation regarding projectable sections of vector bundles. Let 
\[
\xymatrix{
E_1  \ar[d]  \ar[r]^\psi  &  E_2 \ar[d] \\
M_1  \ar[r]_f  &  M_2
}
\]
be a morphism of vector bundles. Denote by 
\begin{align*}
\bar \psi : \Gamma(E_1) & \to \Gamma(f^\ast E_2) \\ 
f^\ast : \Gamma(E_2) & \to \Gamma(f^\ast E_2)
\end{align*} 
the canonical linear maps determined by $\psi$. We denote by $\Gamma(E_1)^\psi$ the vector space of $\psi$\textbf{-projectable sections} of $E_1$, i.e.
\[
\Gamma(E_1)^\psi \equiv \left\{ \xi \in \Gamma(E_1) \; | \; \text{there exists} \; \xi' \in \Gamma(E_2) \; \text{such that} \; \psi \xi = \xi'  f\right\}
\]
We have:
\[
\xi \text{ is projectable } \Leftrightarrow \;
\bar\psi(\xi) \in \mathrm{Im}\left(f^\ast\right) \; \Leftrightarrow \;
\psi(\xi(x)) = \psi(\xi(y)) \text{ whenever } f(x)=f(y)
 \]
Assume now that $f$ is \textbf{surjective}.
If $\xi \in \Gamma(E_1)$ then the surjectivity of $f$ implies that there exists at most one section $\xi' \in \Gamma(E_2)$ such that $\psi \xi = \xi'f$. It follows that there is a natural linear map
\[
\psi_\ast : \Gamma(E_1)^\psi \to \Gamma(E_2)
\]
We denote by $\Gamma(E_1)^\psi \hookrightarrow \Gamma(E_2)$ the inclusion map. The following is immediate:
\begin{lemma}
\label{lemma: fibre product of projectable sections}
The following diagram commutes and is a fibre product of vector spaces:
\[
\xymatrix{
\Gamma(E_1)^\psi  \ar@{^{(}->}[d]  \ar[r]^{\psi_\ast}   &   \Gamma(E_2)  \ar[d]^{f^\ast} \\
\Gamma(E_1)  \ar[r]_-{\bar \psi}  &   \Gamma(f^\ast E_2)
}
\]
\end{lemma}

\subsection{Projectable sections of $\mathcal{VB}$-groupoids}

We will prove an analogue of Lemma \ref{lemma: fibre product of projectable sections} for morphisms of $\mathcal{VB}$-groupoids. We continue in the setting defined in Section \ref{subsection: VB notation} above, but assume from now on that $\Phi$ is \textbf{surjective on objects, morphisms, and on composable pairs}. This implies that the same is true of $\phi$. 

\begin{definition}
We define the vector space $\Gamma_{mult}(\mathcal V)^\Phi$ of $\Phi$-projectable multiplicative sections of $\mathcal V$ to be the intersection
\[
\Gamma_{mult} (\mathcal V)^\Phi \equiv \Gamma_{mult}(\mathcal V) \cap \Gamma(\mathcal V)^\Phi
\]
\end{definition}

\begin{proposition}
\label{prop: projecting mult sections}
 Let $V \in \Gamma_{mult}(\mathcal V)^\Phi$ and $v \in \Gamma(E)$ be the object component of $V$. Then:
 \begin{enumerate}
  \item  $\Phi_\ast V \in \Gamma_m(\mathcal W)$.
 \item $v \in \Gamma(E)^{\Phi_0}$.
\item $(\Phi_0)_\ast v$ is the object component of $\Phi_\ast V$.
 \end{enumerate}
\end{proposition}

\begin{proof}
1. Let $(h,h') \in \mathcal H_2$ be a composable pair of morphism in $\mathcal H$. By assumption $\phi$ is surjective on composable pairs, so we can choose $(g,g') \in \mathcal G_2$ such that $(\phi(g),\phi(g')) = (h,h')$. This implies that $\phi(gg') = hh'$. Then
\begin{align*}
(\Phi_\ast V) (h) \cdot (\Phi_\ast V)(h') 
& =  \left(\Phi \circ V(g) \right)  \cdot \left( \Phi \circ V(g')\right) \\
& = \Phi   \left(V(g) \cdot V(g') \right) \\
& = \Phi \left( V(gg')\right) \\
& = \left(\Phi_\ast V \right)(hh')
\end{align*}
which shows that $\Phi_\ast V$ is multiplicative.

2.  If $x,y \in M$ are such that $\phi_0(x)=\phi_0(y)$ then $\phi(1_x) = \phi(1_y)$ and therefore $V(1_x) = V(1_y)$ by part 1. Then
\begin{align*}
v (x)  &  = v s(1_x) \\
& = \tilde s v (1_x) \\
& = \tilde s v(1_y) \\
& = v s (1_y) \\
& = v (y)
\end{align*}
which shows that $v \in \Gamma(E)^{\Phi_0}$. 

3. Since we have shown that $\Phi_\ast V$ is multiplicative it is sufficient to show that 
\[
\tilde s' \circ \Phi_\ast V = (\Phi_0)_\ast v  \circ s'
\]

Let $h \in \mathcal H$. By assumption $\phi$ is surjective, so we can choose $g \in \mathcal G$ such that $\phi(g)=h$. This implies that $\phi_0(s(g)) = s'(h)$. Then:
\begin{align*}
\left(\tilde s' \circ \Phi_\ast  V   \right) (h)  &  =  \tilde s'  \Phi \left( V(g) \right) \\
& = \Phi_0 \tilde s V(g) \\
& = \Phi_0 v s (g) \\
& = \left( (\Phi_0)_\ast v \right) (s'(h)) \\
& = \left( (\Phi_0)_\ast  v \circ s' \right) (h) \\
\end{align*}
as required.
\end{proof}
 It follows from Proposition \ref{prop: projecting mult sections} that $\Phi_\ast$ restricts to a linear map
 \[
 \Phi_\ast : \Gamma_m(\mathcal V)^\Phi \to \Gamma_m(\mathcal W)
 \]
 The following Lemma is an analogue of Lemma \ref{lemma: fibre product of projectable sections} above:

 \begin{lemma}
The following diagram commutes and is a fibre product of vector spaces:
\begin{equation}
\label{diag: fibre product of mult sections}
\xymatrix{
\Gamma_{mult}(\mathcal V)^\Phi  \ar@{^{(}->}[d]  \ar[r]^{\Phi_\ast}   &   \Gamma_{mult}(\mathcal W)  \ar[d]^{\phi^\ast} \\
\Gamma_{mult}(\mathcal V)  \ar[r]_-{\bar \Phi}  &   \Gamma_{mult}(\phi^\ast \mathcal W)
}
\end{equation}
\end{lemma}

\begin{proof}
It follows from Lemma \ref{lemma: fibre product of projectable sections} that (\ref{diag: fibre product of mult sections}) commutes. The fact that (\ref{lemma: fibre product of projectable sections}) is a fibre product holds for essentially the same reason as Lemma \ref{lemma: fibre product of projectable sections}: as $\phi$ is surjective the map $\phi^\ast$ in (\ref{diag: fibre product of mult sections}) is injective, from which it follows that the fibre product of the linear maps $\bar \Phi$ and $\phi^\ast$ is isomorphic to the preimage $\bar \Phi^{-1} \left( \mathrm{Im} \left(\phi^\ast\right) \right)$:
\begin{align*}
\Gamma_{mult}\left(\mathcal V\right) \times_{\Gamma_{mult}\left(\phi^\ast \mathcal W\right)}  \Gamma_{mult}\left(\mathcal W\right)
& =  \bar \Phi^{-1} \left( \mathrm{Im} \left(\phi^\ast\right) \right) \\
& = \Gamma_{mult}\left(\mathcal V\right) \cap \Gamma\left(\mathcal V\right)^\Phi \\
& = \Gamma_{mult}\left(\mathcal V\right)^\Phi
\end{align*}
Moreover, the two projections out of the fibre product can then be identified with the inclusion $\Gamma_{mult}\left(\mathcal V\right)^\Phi \hookrightarrow \Gamma_{mult}\left(\mathcal V\right)$ and the map $\Phi_\ast.$
\end{proof}

Applying Lemma \ref{lemma: fibre product of projectable sections} to the morphism $\Phi_C : C \to D$ gives

\begin{lemma}
The following diagram is a fibre product of vector spaces:
\begin{equation}
\label{diag: fibre product of core sections}
\xymatrix{
\Gamma(C)^{\Phi_C}  \ar@{^{(}->}[d]  \ar[r]^{(\Phi_C)_\ast}   &   \Gamma(D)  \ar[d]^{\phi_0^\ast} \\
\Gamma(C)  \ar[r]_-{\bar \Phi_C}  &   \Gamma(\phi_0^\ast F)
}
\end{equation}
\end{lemma}

\subsection{The projectable subcomplex}
 
 We will show that $\Gamma(C)^{\Phi_C}$ and $\Gamma_{mult}(\mathcal V)^\Phi$ form a subcomplex $C_{mult}^\bullet(\mathcal V)^\Phi$ of the complex $C_{mult}^\bullet(\mathcal V)$ of multiplicative sections of $\mathcal V$.

\begin{proposition}
\label{prop: fibre product construction of complexes}
The fibre product of complexes
\[
\xymatrixcolsep{4pc}
\xymatrix{
C_{mult}^\bullet (\mathcal V)  \times_{C_{mult}^\bullet (\phi^\ast \mathcal W)}
C_{mult}^\bullet ( \mathcal W) 
\ar[d]   \ar[r]  
&  C_{mult}^\bullet ( \mathcal W) \ar[d]^{\phi^\ast}  \\
C_{mult}^\bullet (\mathcal V)  \ar[r]_{\bar \Phi} 
& C_{mult}^\bullet (\phi^\ast \mathcal W)
}
\]
is isomorphic to the complex
\[
\delta : \Gamma(C)^{\Phi_C}  \to \Gamma_{mult}(\mathcal V)^\Phi
\]
defined by the restriction of $\delta$ to $\Gamma(C)^{\Phi_C}$. In particular, if $c \in \Gamma(C)^{\Phi_C}$ then $\delta(c) \in \Gamma_{mult}(\mathcal V)^\Phi$.
\end{proposition}

\begin{proof}
Consider the following diagram:
\[
\xymatrixcolsep{3pc}
\xymatrix{
\Gamma(C) \times_{\Gamma(f^\ast D)} \Gamma(D)  \ar[r]^{(\delta_\mathcal V,\delta_\mathcal W)}  \ar[d]_{\mathrm{pr}_1}  &
\Gamma(\mathcal V)  \times_{\Gamma(\phi^\ast \mathcal W)} \Gamma_m(\mathcal W)  \ar[d]^{\mathrm{pr}_1}  \\
\Gamma(C)^{\Phi_C}  \ar[r]_{\delta_\mathcal V}  &  \Gamma_{mult}(\mathcal V)^\Phi
}
\]
The projections $\mathrm{pr}_1$ take values in $\Gamma\left(C\right)^{\Phi_C}$ and $\Gamma_{mult}\left(\mathcal V\right)^\Phi$ and are isomorphisms onto those subspaces by Lemmas \ref{diag: fibre product of mult sections} and \ref{diag: fibre product of core sections}. It is immediate that the diagram commutes.
\end{proof}

 \begin{definition}
The complex 
\[
C_{mult}^\bullet (\mathcal V)^\Phi = \left(\delta : \Gamma(C)^{\Phi_C}  \to \Gamma_{mult}(\mathcal V)^\Phi\right)
\] 
is the \textbf{subcomplex of projectable sections}.
\end{definition}

\begin{theorem}
\label{thm: fibre product of qisoms}
Let $\Phi: \mathcal V \to \mathcal W$ be a $\mathcal{VB}$-Morita map covering a morphism $\phi:\mathcal G \to \mathcal H$ and assume that $\Phi$ is surjective on objects. Then the following diagram is a fibre product of 2-term complexes in which every morphism is a quasi-isomorphism:
\begin{equation}
\label{diagram: fibre product of complexes}
\xymatrixcolsep{4pc}
\xymatrix{
C_{mult}^\bullet (\mathcal V)^\Phi
\ar@{^(->}[d]   \ar[r]^{\Phi_\ast}  
&  C_{mult}^\bullet ( \mathcal W) \ar[d]^{\phi^\ast}  \\
C_{mult}^\bullet (\mathcal V)  \ar[r]_{\bar \Phi} 
& C_{mult}^\bullet (\phi^\ast \mathcal W)
}
\end{equation}
\end{theorem}

Theorem \ref{thm: fibre product of qisoms} follows from the following two Lemmas, the proofs of which are straightforward:
\begin{lemma}
\label{lemma: surjectivity in degree 1}
Let $f:V^\bullet \to W^\bullet$ be a quasi-isomorphism of 2-term complexes. If $f^0$ is surjective then so is $f^1$.
\end{lemma}

\begin{lemma}
\label{lemma: pullback of surjective quasi-isomorphisms}
Let 
\[
\xymatrix{
U^\bullet \ar[d] \ar[r]^{f'} &   V^\bullet \ar[d] \\
W^\bullet \ar[r]_f  &  Z^\bullet
}
\]
be a fibre product of cochain complexes. If $f$ is a surjective quasi-isomorphism then so is $f'$.
\end{lemma}

\begin{proof}\emph{(of Theorem \ref{thm: fibre product of qisoms})}
The fact that $\Phi$ is a Morita map and surjective on objects implies that it is also surjective on morphisms and composable pairs. It then follows from Proposition  \ref{prop: fibre product construction of complexes} that (\ref{diagram: fibre product of complexes}) is a fibre product of 2-term complexes.  By Proposition \ref{prop: VB qisoms} the morphisms $\bar \Phi$ and $\phi^\ast$ are quasi-isomorphisms. The degree 0 component of $\bar \Phi$ is the map
\[
\bar \Phi_C : \Gamma(C) \to \Gamma\left(\phi_0^\ast D\right)
\]
which is surjective because $\bar \Phi_C$ and $\phi_0$ are surjective. It then follows from Lemma \ref{lemma: surjectivity in degree 1} that the degree 1 component of $\bar \Phi$ is surjective, and then by Lemma \ref{lemma: pullback of surjective quasi-isomorphisms} that the map $\Phi_\ast$ is a quasi-isomorphism. As three of the four morphisms in the commutative square (\ref{diagram: fibre product of complexes}) are then quasi-isomorphisms, it follows that the same is true of the inclusion $C^\bullet_{mult}\left(\mathcal V\right)^\Phi \hookrightarrow C^\bullet_{mult}\left(\mathcal V\right)$.
\end{proof}

\section{Sections of $\mathcal{LA}$-groupoids}
\label{section: sections of LA}

\subsection{Lie 2-algebra structure}

Let $\mathcal{V}$ be an $\la$-groupoid over $\G\rr M$ with side Lie algebroid $E\rmap M$ and core $C$. In this section we will see that the complex of multiplicative sections $\delta:\Gamma(C)\rmap \Gamma_{mult}(\mathcal{V})$ is part of a crossed module of Lie algebras. In particular, the category $\mathrm{Sec}(\G,\mathcal{V})$ of multiplicative sections of an $\la$-groupoid inherits the structure of a Lie 2-algebra.

Recall that, as explained in subsection \ref{LAgroupoids}, the core bundle has a canonical structure of Lie algebroid. Also, for any $\la$-groupoid $\mathcal{V}$ the space of multiplicative sections $\Gamma_{mult}(\mathcal{V})$ is closed under the Lie bracket $[\cdot,\cdot]_{\mathcal{V}}$ on $\Gamma(\mathcal{V})$. Hence, we have a pair of Lie algebras $(\Gamma(C),\Gamma_{mult}(\mathcal{V}))$ naturally associated to the $\la$-groupoid $\mathcal{V}$.

The following result due to Mackenzie (c.f. Lemma 2.12 in \cite{mackII}) will be used to prove Theorem \ref{thm:crossed} below.  Recall that a \textbf{star section} of $\mathcal{V}$ is a pair $(V,e)$ where $V\in\Gamma(\mathcal{V})$ and $e\in\Gamma(E)$ satisfy $\tilde{s}\circ V=e\circ s$ and $V\circ 1=\tilde{1}\circ e$. For instance, every multiplicative section is a star section. As shown in \cite{mackII}, given a star section $(V,e)$ and a section $c\in \Gamma(C)$, define 

\begin{equation}\label{derivationmackII}
D_V(c):=[V,c^r]_{\mathcal{V}}\circ 1\in \Gamma(C)
\end{equation}

\begin{lemma}\label{lemmamackII}(\cite{mackII}) Let $\mathcal{V}\rr E$ be an $\la$-groupoid over $\G\rr M$. Given star sections $(V,e)$ and $(V',e')$ and a section $c\in\Gamma(C)$, the following hold:

\begin{itemize}

\item[i)] If $W\in \Gamma(\mathcal{V})$ is any section with $W\circ 1=c$, then $D_V(c)=[V,W]\circ 1$.

\item[ii)] $D_{[V,V']}(c)=D_{V}(D_{V'}(c))-D_{V'}(D_{V}(c)).$

\item[iii)] For any $f\in C^{\infty}(M)$, $D_V(fc)=fD_V(c)+(\mathcal{L}_{\rho_E(e)}f)c.$

\end{itemize} 

\end{lemma}

We will use the previous lemma in the case of multiplicative sections, concluding that the Lie algebra $\Gamma_{mult}(\mathcal{V})$ acts on the Lie algebra $\Gamma(C)$ by derivations. Any multiplicative section $V:\G\rmap \mathcal{V}$ covering $e:M\rmap E$ defines a derivation $D_V:\Gamma(C)\rmap \Gamma(C)$ with symbol $\rho_E(e)$.

\begin{theorem}\label{thm:crossed}

The pair of Lie algebras $(\Gamma(C),\Gamma_{mult}(\mathcal{V}))$ has a natural crossed module structure given by

\begin{equation}
\xymatrix{\Gamma(C)\ar[r]^{\delta}&\Gamma_{mult}(\mathcal{V})\ar[r]^{D}&\mathrm{Der}(\Gamma(C))},
\end{equation}
where $\delta(c)=c^r-c^l$ and $D:\Gamma_{mult}(\mathcal{V})\rmap \mathrm{Der}(\Gamma(C))$ is defined by \eqref{derivationmackII}.

\end{theorem}

\begin{proof}

Note that ii) in Lemma \ref{lemmamackII} says that $D:\Gamma_{mult}(\mathcal{V})\rmap \mathrm{Der}(\Gamma(C))$ is a Lie algebra morphism. Now we prove that $\delta:\Gamma(C)\rmap \Gamma_{mult}(\mathcal{V})$ is also a Lie algebra morphism. In fact

\begin{align*}
\delta([c_1,c_2])=&[c_1,c_2]^r-[c_1,c_2]^l\\
=&[c_1^r,c_2^r]+[c_1^l,c_2^l]\\
=&[c_1^r-c_1^l,c_2^r-c_2^l]\\
=& [\delta(c_1),\delta(c_2)],
\end{align*}

\noindent where we have used the facts that $[c_1,c_2]^l=-[c^l_1,c^l_2]$, which follows from the fact that both the multiplication and inversion of $\mathcal{V}\rr E$ are Lie algebroid morphisms, and $[c_1^r,c_2^l] = [c_1^l,c_2^r]=0$, which follows from the definition of $c^r$ and $b^l$ together with the fact that $\tilde{m}:\mathcal{V}^{(2)}\rmap \mathcal{V}$ is a Lie algebroid morphism. Let us prove now the first condition of a crossed module, that is $D_{\delta(c)}=\mathrm{ad}_{c}$ for every $c\in\Gamma(C)$. For that, take a section $b\in\Gamma(C)$, then

\begin{align*}
[c^r-c^l,b^r]_{\mathcal{V}}=&[c^r,b^r]_{\mathcal{V}}-[c^l,b^r]_{\mathcal{V}}\\
=& [c,b]^r-[c^l,b^r]_{\mathcal{V}}\\
=&[c,b]^r.
\end{align*}

\noindent In the last equality we used that $[c^r,b^l]=0$.  Therefore, we conclude that $[\delta(c),b^r]=[c,b]^r$, showing that $D_{\delta(c)}=\mathrm{ad}_{c}$. We proceed now to prove the second condition of a crossed module, namely $\delta\circ D_V=\mathrm{ad}_V\circ \delta$ for every multiplicative section $V\in\Gamma_{mult}(\mathcal{V})$. On the one hand, for every section $c\in \Gamma(C)$ we have 

\begin{align*}
\delta(D_V(c))=&(D_V(c))^r-(D_V(c))^l\\
=& [V,c^r]_{\mathcal{V}}-(D_V(c))^l.
\end{align*}

\noindent On the other hand, we have

\begin{align*}
[V,\delta(c)]_{\mathcal{V}}=[V,c^r-c^l]_{\mathcal{V}}=[V,c^r]_{\mathcal{V}}-[V,c^l]_{\mathcal{V}}.
\end{align*}

\noindent We conclude that $\delta(D_V(c))=[V,\delta(c)]$ if and only if $(D_V(c))^l=[V,c^l]_{\mathcal{V}}$. The last identity follows from applying the inversion $\tilde{i}:\mathcal{V}\rmap \mathcal{V}$ to both sides of the identity $(D_Vc)^r=[V,c^r]_{\mathcal{V}}$.
\end{proof}

We have shown in Theorem \ref{2vect} that the category of multiplicative sections $\mathrm{Sec}(\G,\mathcal{V})$ of a $\vb$-groupoid $\mathcal{V}\rr E$ over $\G\rr M$, has a natural structure of 2-vector space. In the special case of an $\la$-groupoid the category $\mathrm{Sec}(\G,\mathcal{V})$ inherits a Lie 2-algebra structure. 

\begin{theorem}\label{Lie2}

Let $\mathcal{V}\rr E$ be an $\la$-groupoid over $\G\rr M$. The 2-vector space $\Gamma(C)\oplus \Gamma_{mult}(\mathcal{V})\rr \Gamma_{mult}(\mathcal{V})$ has a Lie 2-algebra structure induced by the crossed module of Theorem \ref{thm:crossed}.

\end{theorem} 

The following result is immediate.

\begin{corollary}\label{cohomologyliealgebra}

The cohomology groups $H^{0}(C^\bullet_{mult}(\mathcal{V}))$ and $H^{1}(C^\bullet_{mult}(\mathcal{V}))$ carry natural Lie algebra structures induced by the crossed module of Lie algebras $(\Gamma(C),\Gamma_{mult}(\mathcal{V}))$.

\end{corollary}

\subsection{Examples}

Let us see some examples.

\begin{example}(Lie 2-algebras)
We saw in Example \ref{example: sections of 2-vector space} that the category of multiplicative sections of a 2-vector space $V\rr W$ is $V\rr W$ itself. One can easily see that if $V\rr W$ is a Lie 2-algebra, then the crossed module of Theorem \label{thm:crossed} is the crossed module associated to $V_1 \rightrightarrows V_0$, and the Lie 2-algebra of Theorem \ref{Lie2} coincides with $V\rr W$.
\end{example}

\begin{example}

Let $\G\rr M$ be a Lie groupoid with Lie algebroid $A$. Theorem \ref{Lie2} implies that the complex of multiplicative sections of the tangent $\la$-groupoid $T\G\rr TM$ is $\delta:\Gamma(A)\rmap \mathfrak{X}_{mult}(\G)$ where $\delta(a)=a^r-a^l$ is part of a crossed module of Lie algebras $(\Gamma(A),\mathfrak{X}_{mult}(\G))$. As a consequence, the category of multiplicative vector fields $\Gamma(A)\oplus \mathfrak{X}_{mult}(\G)\rr \mathfrak{X}_{mult}(\G)$ inherits a Lie 2-algebra structure, given by

$$\Cour{(a,X),(b,Y)}=([a,b]_A+D_X(b)-D_Y(a),[X,Y]),$$

\noindent where $D:\mathfrak{X}_{mult}(\G)\rmap \mathrm{Der}(\Gamma(A))$ is given as in \eqref{derivationmackII}.
\end{example}

\begin{remark}
The category of multiplicative vector fields $\mathrm{Vect}(\G)$ on a Lie groupoid was introduced by Hepworth in \cite{hepworth}. It was conjectured in \cite{hepworth} that the category $\mathrm{Vect}(\G)$ has a Lie 2-algebra structure. The example above says that Hepworth's conjecture is true. However, Hepworth's conjecture is a consequence of a much more general result, namely Theorem \ref{Lie2}, which states that the category of multiplicative sections of \emph{any $\la$-groupoid} has the structure of a Lie 2-algebra. 
\end{remark}

\begin{example}(Poisson groupoids)

Let $(\G,\pi_\G)$ be a Poisson groupoid. As observed in Example \ref{multiplicativeforms}, the complex of multiplicative sections in this case is $\delta:\Omega^1(M)\rmap \Omega^1_{mult}(\G); \alpha\mapsto t^*\alpha-s^*\alpha$. It is well known that the base $M$ of any Poisson groupoid $(\G,\pi_\G)$ inherits a unique Poisson structure $\pi_M$ making the source map into a Poisson map. As a consequence, $\Omega^1(M)$ has a Lie algebra structure which coincides with the one induced on the space of sections of the core of the cotangent $\la$-groupoid. As shown in \cite{mackII}, for any multiplicative 1-form $\omega\in\Omega^1_{mult}(\G)$ covering $\xi\in \Gamma(A^*)$, the covariant differential operator $D_\omega:\Omega^1(M)\rmap \Omega^1(M); \alpha\mapsto [\omega,\alpha^r]_{\pi_\G}$ (defined generally by \eqref{derivationmackII}), has symbol $\rho_{A^*}(\xi)$ where $\rho_{A^*}:A^*\rmap TM$ denotes the anchor map of the dual Lie algebroid $A^*$. Since $\alpha^r=t^*\alpha$, as explained in Example 2.18 of \cite{mackII}, the following holds

$$[\omega,t^*\alpha]_\G=t^*(D_\omega \alpha).$$

\noindent Therefore, the Lie 2-algebra of multiplicative 1-forms on $\G$ is given by $\Omega^1(M)\oplus \Omega^1_{mult}(\G)\rr \Omega^1_{mult}(\G)$ with the following bracket

$$\Cour{(\alpha,\omega),(\beta,\theta)}=([\alpha,\beta]_{\pi_M}+D_\omega \beta - D_\theta \alpha,[\omega,\theta]_{\pi_{\G}}).$$

\end{example}

\subsection{Dual Lie brackets}

Assume now that $\mathcal{V}\rr E$ is an $\la$-groupoid over $\G\rr M$. The dual vector bundle $\mathcal{V}^*\rmap \G$ has a linear Poisson structure which is also multiplicative with respect to the groupoid structure $\mathcal{V}^*\rr C^*$, making $\mathcal{V}^*$ into a \emph{PVB-groupoid}. Here PVB stands for \emph{Poisson vector bundle}. Since $\mathcal{V}^*$ comes equipped with a linear Poisson bracket $\{\cdot,\cdot\}_{\mathcal{V}^*}$, the space of linear 1-cochains $C^{\infty}_{lin}(\mathcal{V}^*)$ inherits a Lie algebra structure with Lie bracket $\{\cdot,\cdot\}_{\mathcal{V}^*}$. Additionally, the fact that $\{l_V,l_W\}_{\mathcal{V}^*}=l_{[V,W]_{\mathcal{V}}}$, implies that the space of linear 1-cocycles $Z^1_{lin}(\mathcal{V}^*)\subseteq C^{\infty}_{lin}(\mathcal{V}^*)$ is a Lie subalgebra. As a result we conclude that the isomorphism of Proposition \ref{mult1cocycles} 

$$\Gamma_{mult}(\mathcal{V})\rmap Z^1_{lin}(\mathcal{V}^*),$$

\noindent is a Lie algebra isomorphism. 

Since $\mathcal{V}\rr E$ is an $\la$-groupoid, the linear Poisson groupoid $\mathcal{V}^*\rr C^*$ induces a Poisson structure on the base $C^*$ which is also linear. Indeed, this is the linear Poisson structure dual to the Lie algebroid structure on the core $C\rmap M$. In particular, the canonical identification $\Gamma(C)\rmap C^{\infty}_{lin}(C^*)$ is a Lie algebra isomorphism. As a result we get.

\begin{proposition}

The isomorphism of 2-term complexes between $\delta:\Gamma(C)\rmap \Gamma_{mult}(\mathcal{V})$ and $d:C^{\infty}_{lin}(C^*)\rmap Z^1_{lin}(\mathcal{V}^*)$ is given by Lie algebra isomorphisms at each degree. 
\end{proposition}

\subsection{$\mathcal{LA}$-Morita maps}
\label{section: Morita equivalence LA}

\begin{definition}

Let $\mathcal{V}_1$ and $\mathcal{V}_2$ be $\mathcal{LA}$-groupoids over $\G_1$ and $\G_2$, respectively. A morphism of $\mathcal{LA}$-groupoids $\left(\Phi,\phi\right):\mathcal V_1 \to \mathcal V_2$ is called an \textbf{$\mathcal{LA}$-Morita map} if $\Phi$ is a weak equivalence of Lie groupoids.
\end{definition}

\begin{example}
\label{example: T is an LA morita map}
Let $\phi:\mathcal G \to \mathcal G'$ be a Morita map. Then the differential $T\phi: T\mathcal G \to T\mathcal G'$ is an $\mathcal{LA}$-Morita map. (This follows from the fact that the tangent functor $T$ preserves transverse fibre products and surjective submersions, or from the results of \cite{dHO}.)
\end{example}

\begin{definition}
We say that $\mathcal{LA}$-groupoids $\mathcal V$ and $\mathcal V'$ are \emph{Morita equivalent} if there exists a third $\mathcal{LA}$-groupoid $\mathcal W$ and $\mathcal{LA}$-Morita morphisms $\Phi:\mathcal W \to \mathcal V$ and $\Phi':\mathcal W \to \mathcal V'$ as in the following diagram:
\[
\xymatrix{
& \mathcal W \ar[dl]_{\Phi'}  \ar[dr]^{\Phi'} & \\
\mathcal V  & & \mathcal V'
}
\]
\end{definition}

\begin{proposition} 
Morita equivalence for $\mathcal{LA}$-groupoids is an equivalence relation.
\end{proposition}

\begin{proof}
Using Proposition \ref{prop: existence of LA weak fibre products} the proof is identical to the analogous statement for Lie groupoids (\cite{MoerdijkMrcun}).
\end{proof}

\begin{proposition} 
\label{prop: surjective LA Morita maps}
If $\mathcal V$ and $\mathcal V'$ are Morita equivalent $\mathcal{LA}$-groupoids then there exists an $\mathcal{LA}$-groupoid $\mathcal W$ and $\mathcal{LA}$-Morita maps $\mathcal W \to \mathcal V$ and $\mathcal W \to \mathcal V'$ which are surjective on objects, morphisms and composable pairs.
\end{proposition}

\begin{proof}
Using Proposition \ref{prop: existence of LA weak fibre products} the proof is identical to the analogous statement for Lie groupoids (\cite{MoerdijkMrcun}).
\end{proof}

\subsection{Morita invariance of $C^\bullet_{mult}(\mathcal V)$ for $\mathcal{LA}$-groupoids}

If $\mathcal V$ and $\mathcal V'$ are Morita equivalent $\mathcal{LA}$-groupoids then they are in particular Morita equivalent as $\mathcal{VB}$-groupoids, and so Theorem \ref{moritainvarianceVB2} shows that $H^\bullet\left(C_{mult}^\bullet(\mathcal V)\right) \cong H^\bullet\left(C_{mult}^\bullet(\mathcal V)\right)$ as graded vector spaces. However, this isomorphism will not in general preserve the graded Lie brackets. This issue arises from the fact that even if $\mathcal W$ is a Lie algebroid, the $\mathcal{VB}$-groupoid $\phi^\ast \mathcal W^\ast$ in (\ref{VBmappullback2})  does not carry any natural Lie structure. 

Instead, we prove a result at the level of the complexes. We show that the projectable subcomplex $C^\bullet_{mult}(\mathcal V)^\Phi$ is a sub crossed module of $C^\bullet_{mult}(\mathcal V)$, and that the inclusion and projection morphisms in (\ref{diagram: fibre product of complexes}) are quasi-isomorphisms of crossed modules. We use this to show that if $\mathcal V$ and $\mathcal V'$ are Morita equivalent $\mathcal{VB}$-groupoids then the crossed modules $C^\bullet_{mult}\left(\mathcal V\right)$ and $C^\bullet_{mult}\left(\mathcal V'\right)$ are isomorphic in the derived category of crossed modules.

\begin{proposition}
\label{prop: projection and inclusion are morphisms}
Let $\Phi: \mathcal V \to \mathcal W$ be an $\mathcal{LA}$-Morita map and assume that $\Phi$ is surjective on objects. Then:
\begin{enumerate}
\item The projectable subcomplex $C^\bullet_{mult}\left(\mathcal V\right)^\Phi$ is a sub crossed module of $C^\bullet_{mult}\left(\mathcal V\right)$.
\item The inclusion $C_{mult}^\bullet (\mathcal V)^\Phi \hookrightarrow C_{mult}^\bullet(\mathcal V)$ is a morphism of crossed modules.
\item The projection $\Phi_\ast : C_{mult}^\bullet (\mathcal V)^\Phi \to C_{mult}^\bullet(\mathcal W)$ is a morphism of crossed modules.
\end{enumerate}
\end{proposition}

\begin{proof}
1. Recall from \cite{mackbook} that if $\psi:A \to B$  is a morphism of Lie algebroids, $a,a' \in \Gamma(A)$, $b,b' \in \Gamma(B)$, and $a$ respectively $a'$ projects to $b$ respectively $b'$, then $\left[a,a'\right]$ projects to $\left[b,b'\right]$. As $\Phi:\mathcal V \to \mathcal W$ and $\Phi_C: C \to D$ are morphisms of Lie algebroids it then follows that $\Gamma\left(C\right)^{\Phi_C} \subset \Gamma\left(C\right)$ and $\Gamma_{mult}\left(\mathcal V\right)^\Phi \subset \Gamma_{mult}\left(\mathcal V\right)$ are Lie subalgebras. 

It remains to show that $\Gamma\left(C\right)$ is closed under the action of $\Gamma_{mult}\left(\mathcal V\right)$. Let $V \in \Gamma_{mult}\left(\mathcal V\right)^\Phi$, $\Phi_\ast V=W$, $c \in \Gamma\left(C\right)^{\Phi_C}$, and $\left(\Phi_C\right)_\ast c = d$. As shown in the proof of Proposition \ref{prop: pullback chain map}, $\Phi_\ast c^r = d^r$. Therefore $\left[ V,c^r\right]$ projects to $\left[W,d^r\right]$, which shows that $D_V(c)$ projects to $D_W(d)$. In particular, $D_V(c)$ is projectable.

2.  This follows immediately from 1.

3. This follows from the proof of 1: the projection maps $\left(\Phi_C\right)_\ast$ and $\Phi_\ast$ are Lie algebra morphisms, and 
\begin{equation*}
\left(\Phi_C\right)_\ast \left(D_V (c)\right) = D_{\Phi_\ast V} \left( \Phi_C\right)_\ast c
\end{equation*}
\end{proof}

\begin{theorem}
\label{theorem: quis from surjective maps}
If $\Phi: \mathcal V \to \mathcal W$ is an $\mathcal{LA}$-Morita map that is surjective on objects then the maps 
\[
\xymatrixcolsep{3pc}
\xymatrix{
& C_{mult}^\bullet(\mathcal V)^\Phi \ar@{_(->}[dl] \ar[dr]^{\Phi_\ast} \\
 C_{mult}^\bullet(\mathcal V) && C_{mult}^\bullet(\mathcal W)
 }
 \]
are quasi-isomorphisms of crossed modules.
\end{theorem}

\begin{proof}
By Proposition \ref{prop: projection and inclusion are morphisms} the inclusion and the projection $\Phi_\ast$ are morphisms of crossed modules, and by Theorem \ref{thm: fibre product of qisoms} they are both quasi-isomorphisms.
\end{proof}

\begin{theorem} 
\label{theorem: invariance of Lie 2-algebras}
If $\mathcal V$ and $\mathcal V'$ are Morita equivalent $\mathcal{LA}$-groupoids then the crossed modules $C_{mult}^\bullet(\mathcal V)$ and $C_{mult}^\bullet (\mathcal V')$ are isomorphic in the derived category of crossed modules.
\end{theorem}

\begin{proof}
By Proposition \ref{prop: surjective LA Morita maps} there exists an $\mathcal{LA}$-groupoid $\mathcal W$ and $\mathcal{LA}$-Morita maps $\Phi:\mathcal W \to \mathcal V$ and $\Psi:\mathcal W \to \mathcal V'$, both of which are surjective on objects. By Theorem \ref{theorem: quis from surjective maps} this induces a zig-zag of quasi-isomorphisms of crossed modules:
\begin{equation}
\label{eqn: zigzag of qisoms}
\xymatrix{
C_{mult}^\bullet(\mathcal V)  & 
\ar[l]_{\Phi_\ast} C_{mult}^\bullet(\mathcal W)^\Phi \ar@{^(->}[r] & 
C_{mult}^\bullet(\mathcal W)   & 
\ar@{_(->}[l] C_{mult}^\bullet(\mathcal W)^\Psi \ar[r]^{\Psi_\ast}  & 
C_{mult}^\bullet(\mathcal V')  & 
}
\end{equation}
\end{proof}

\begin{corollary}\label{cor:invariancecohomologyliealgebra}
If $\mathcal V$ and $\mathcal V'$ are Morita equivalent $\mathcal{LA}$-groupoids then there are Lie algebra isomorphisms $H^0(C_{mult}(\mathcal{V}))\cong H^0(C_{mult}(\mathcal{V}'))$ and $H^1(C_{mult}(\mathcal{V}))\cong H^1(C_{mult}(\mathcal{V}'))$. In particular, an $\mathcal{LA}$-Morita map $\Phi:\mathcal{V} \to \mathcal{V}'$ that is surjective on objects determines such a pair of isomorphisms. 
\end{corollary}

\subsection{$L_\infty$-morphisms}
\label{subsection: Linfinity}

If we consider crossed modules as differential graded Lie algebras concentrated in degrees $0$ and $1$, then we can speak of $L_\infty$-morphisms between them. See \cite{Noohi} for this perspective on crossed modules, and, for example, \cite{Kontsevich} for the notions of $L_\infty$-morphism.

\begin{theorem}
\label{theorem: L infinity map}
If $\mathcal V$ and $\mathcal V'$ are Morita equivalent $\mathcal{LA}$-groupoids then there exists an $L_\infty$ quasi-isomorphism 
\[
C_{mult}^\bullet\left(\mathcal V\right) \xrightarrow{\sim} C_{mult}^\bullet \left(\mathcal V'\right)
\] 
\end{theorem}

\begin{proof}
This follows from the fact that if $f:L_1 \to L_2$ is an $L_\infty$-quasi-isomorphism then there exists an $L_\infty$-quasi-isomorphism $f':L_2 \to L_1$ which, at the level of cohomology, is the inverse of $f_1$ (see, for example, \cite{Kontsevich}, Theorem 4.6). We can apply this result to the left pointing morphisms in (\ref{eqn: zigzag of qisoms}) and then compose the resulting composable chain of $L_\infty$-quasi-isomorphisms, the result being an $L_\infty$-quasi-isomorphism from $C^\bullet_{mult}\left(\mathcal V\right)$ to $C^\bullet_{mult}\left(\mathcal V'\right)$.
\end{proof}


\section{Applications}\label{sec:vectorfields}

In \cite{hepworth} Hepworth introduced vector fields on a differentiable stack $\mathfrak{S}$ as maps $\mathfrak{S}\rmap T\mathfrak{S}$ which are sections of the tangent stack $T\mathfrak{S}$ up to a 2-morphism. Such vector fields form a category $\mathrm{Vect}(\mathfrak{S})$ and one of the main results of \cite{hepworth} establishes that this category is equivalent to the category of multiplicative vector fields on any Lie groupoid presenting $\mathfrak{S}$.  

In this section we apply the results of section \ref{categorymultiplicative} to the special case of the tangent $\la$-groupoid $T\G\rr TM$, proving a conjecture by Hepworth \cite{hepworth} about the Lie 2-algebra structure on the category of vector fields on a differentiable stack. This conjecture was proven independently by Berwick-Evans and Lerman in \cite{BELerman}.

One of the main goals of this section consists of introducing a different notion of vector field on a differentiable stack, which is more geometric than categorical. Several examples are discussed.


\subsection{Multiplicative vector fields on Lie groupoids}

Let $\G\rr M$ be a Lie groupoid. Recall that a multiplicative vector field on $\G$ is by definition a multiplicative section of the tangent $\vb$-groupoid $T\G\rr TM$. In this case, the category of multiplicative sections of $T\G$ coincides with the category of multiplicative vector fields denoted by $\mathrm{Vect}(\G\rr M)$ introduced in \cite{hepworth}.

Given a Lie groupoid $\G\rr M$ with Lie algebroid $A$, the complex of multiplicative vector fields is $\delta:\Gamma(A)\rmap \mathfrak{X}_{mult}(\G);a\mapsto a^r-a^l$, where $a^r$ and $a^l$ denote respectively the right and left invariant vector fields on $\G$ induced by $a\in\Gamma(A)$, see Example \ref{multiplicativefields}.

\begin{remark}
In this section, both the Lie 2-algebra $\mathrm{Vect}(\G\rr M)$ of multiplicative vector fields on a Lie groupoid $\G\rr M$ as well as the complex of multiplicative vector fields $\Gamma(A)\rmap \mathfrak{X}_{mult}(\G)$ will be simply denoted by $L_{\G}$.
\end{remark}

A result of \cite{dHO} shows that a Lie groupoid morphism $\varphi:\G_1\rmap \G_2$ is a Morita map if and ony if the tangent morphism $T\varphi:T\G_1\rmap T\G_2$ is a $\vb$-Morita map. Then, due to Theorem \ref{theorem: invariance of Lie 2-algebras} we have

\begin{proposition}\label{vf:moritainvariance}

If $\mathcal G$ and $\mathcal G'$ are Morita equivalent Lie groupoids then the Lie 2-algebras $L_{\G}$ and $L_{\G'}$ are isomorphic in the derived category of Lie 2-algebras.

\end{proposition}

This allows us to describe the complex of multiplicative vector fields, and hence the Lie 2-algebra of multiplicative vector fields, for concrete examples of Lie groupoids.

\begin{example}(Lie groups)\label{vf:liegroup}

Let $G$ be a Lie group, i.e. a Lie groupoid of the form $G\rr \{*\}$. The tangent groupoid $TG$ is also a Lie group, which is isomorphic to the semidirect product Lie group $G\ltimes \mathfrak{g}$ with respect to the adjoint representation. The core is given by the Lie algebra $\mathfrak{g}$. In this case, a multiplicative vector field can be seen as a smooth map $G\rmap \mathfrak{g}$ which is a 1-cocycle. Hence, the complex of multiplicative vector fields coincides with the 1-truncation $\mathfrak{g}\rmap Z^1(G,\mathfrak{g})$ of the group complex of $G$ with coefficients in the adjoint representation, in agreement with Proposition \ref{1truncationisomorphism}.

\end{example}

A Lie groupoid $s,t:\G\rr M$ is called \textbf{transitive} if $(s,t):\G\rmap M\times M$ is a surjective submersion. If $x\in M$, its \textbf{isotropy group} is $\G_x:=s^{-1}(x)\cap t^{-1}(x)$. It is well-known that if $\G$ is transitive, then the inclusion $\G_x\rmap \G$ is a Morita map for every $x\in M$. Hence, Example \ref{vf:liegroup} describes the category of multiplicative vector fields of any transitive Lie groupoid.

\begin{example}(Manifolds)\label{vf:manifolds}

Any smooth manifold $M$ can be viewed as a Lie groupoid by considering the \textbf{unit groupoid} $M\rr M$, on which all the structure maps are the identity. In this case, the tangent groupoid is also a unit groupoid $TM\rr TM$ and the core bundle is $C=0$. A straightforward computation shows that a multiplicative vector field is just a vector field on $M$. The complex of multiplicative vector fields is $0\rmap \mathfrak{X}(M)$ and hence the Lie 2-algebra of vector fields on $M\rr M$ is the unit Lie 2-algebra $\mathfrak{X}(M)\rr \mathfrak{X}(M)$.

\end{example}

\begin{example}(Submersion groupoids)\label{vf:submersions}

Let $\pi:M\rmap N$ be a surjective submersion. The \textbf{submersion groupoid} denoted $M\times_NM\rr M$ has source and target given by $s(x,y)=y$ and $t(x,y)=x$, and all the other structure maps are the obvious one. The tangent groupoid is the submersion groupoid associated to the surjective submersion $T\pi:TM\rmap TN$. The core bundle is given by the vertical bundle $C=\ker(T\pi)$. In this case, a multiplicative vector field is the same that a $\pi$-projectable vector field. Hence, the category of multiplicative vector fields is given by the 2-vector space associated to the 2-term complex $\Gamma(\ker(T\pi))\rmap \mathfrak{X}(M)^{\pi}$, where $\mathfrak{X}(M)^{\pi}\subseteq \mathfrak{X}(M)$ denotes the subspace of $\pi$-projectable vector fields.

\end{example}

Given a surjective submersion $\pi:M\rmap N$, the canonical map $(M\times_{N}M\rr M)\rmap (N\rr N); (x,y)\mapsto \pi(x),$ is a Morita map. Hence, by Proposition \ref{vf:moritainvariance}, the Lie 2-algebra of multiplicative vector fields $\mathrm{Vec}(M\times_{N}M\rr M)$ is equivalent to the Lie 2-algebra $\mathfrak{X}(N)\rr \mathfrak{X}(N)$ defined in Example \ref{vf:manifolds}. As a particular instance of this situation, associated to any open cover $\mathcal{U}=\{U_i\}$ of a manifold $M$ is the \textbf{\v{C}ech groupoid} which is the submersion groupoid defined by $\pi:\coprod_{i}U_i\rmap M$. Hence, the Lie 2-algebra of vector fields on the \v{C}ech groupoid is equivalent to the Lie algebra of vector fields $\mathfrak{X}(M)$ viewed as a Lie 2-algebra. In other words, a multiplicative vector field on the \v{C}ech groupoid is given by a family of vector fields on each $U_i$ coinciding on the overlaps, which is the same as having a vector field defined on the whole of $M$.

\begin{example}(Transformation groupoid)\label{vf:liegroupaction}

Let $G$ be a Lie group acting smoothly on a manifold $M$. The \textbf{transformation groupoid} is denoted by $G\ltimes M\rr M$ and its arrows are pairs $(g,x)\in G\times M$ with source $s(g,x)=x$ and target $t(g,x)=gx$. A result of Hepworth \cite{hepworth} shows that if $G$ is compact, then the category of multiplicative vector fields on the transformation groupoid $G\ltimes M\rr M$ is equivalent to the 2-vector space associated to the 2-term complex $C^{\infty}(M,\mathfrak{g})^G\rmap \mathfrak{X}(M)^G; (f:M\to \mathfrak{g})\mapsto f(x)_M$, where $C^{\infty}(M,\mathfrak{g})^G=\{f:M\to\mathfrak{g};f(gx)=Ad_g(f(x)), x\in M,g\in G\}$ and $f(x)_{M}$ denotes the infinitesimal generator associated to the Lie algebra element $f(x)\in\mathfrak{g}$ .

\end{example}

Suppose that $G$ is a Lie group which acts freely and properly on a manifold $M$. In this case, the orbit space $M/G$ inherits a unique smooth structure making the canonical projection $\pi:M\rmap M/G$ into a surjective submersion. It is easy to see that the map $\varphi:G\ltimes M\rmap M\times_{M/G} M; (g,x)\mapsto (x,gx)$ is a Lie groupoid isomorphism covering $id:M\rmap M$. Hence, the Lie 2-algebra of multiplicative vector fields on the transformation groupoid $G\ltimes M\rr M$ is \emph{isomorphic} to the Lie 2-algebra of multiplicative vector fields on the submersion groupoid $M\ltimes_{M/G}\rr M$ associated to $\pi:M\rmap M/G$. In particular, if $G$ is a compact Lie group acting freely on $M$, then the Lie 2-algebra of vector fields on $G\ltimes M\rr M$ is equivalent to the Lie algebra of vector fields $\mathfrak{X}(M/G)\cong \mathfrak{X}(M)^G$ seen as a Lie 2-algebra.

\begin{example}(Locally free actions)\label{vf:locallyfree}

Let $G$ be a Lie group acting properly on a manifold $M$. Suppose that the action is locally free, i.e. the isotropy Lie algebra $\mathfrak{g}_x=\{u\in\mathfrak{g}; u_M(x)=0\}$ is trivial, for any $x\in M$. In this case, the transformation groupoid $G\ltimes M\rr M$ has anchor map $\rho:M\times \mathfrak{g}\rmap TM;(x,u)\mapsto u_M(x)$, which is injective and hence $G\ltimes M\rr M$ is a regular groupoid. The tangent groupoid $T(G\ltimes M)\rr TM$ is a regular $\vb$-groupoid with core anchor $\rho:M\times \mathfrak{g}\rmap TM$. Corollary \ref{regularqisotypezero} implies that the complex of multiplicative vector fields on $G\ltimes M$ is quasi-isomorphic to the complex of multiplicative sections of the type zero component of $T(G\ltimes M)\rr TM$, namely $t^*K\oplus s^*\nu\rr \nu$, where $K=\ker(\rho)$ and $\nu=TM/\mathrm{im}(\rho)$. Since the anchor map is injective, we conclude that the complex of multiplicative vector fields on $G\ltimes M\rr M$ is quasi-isomorphic to the complex

$$0\rmap \Gamma_{mult}(s^*\nu).$$

Example \ref{example: Sec of repns} implies that $\Gamma_{mult}(s^*\nu)\cong \Gamma(\nu)^G$, where $\nu\rmap M$ is viewed as a $G$-equivariant vector bundle with respect to the canonical representation of $G$ on $\nu$.

\end{example}

\begin{example}(\'Etale groupoids)\label{vf:etale}

Let $\G\rr M$ be an \'etale Lie groupoid. The tangent groupoid $T\G\rr TM$ has core $C=0$ and hence $T\G$ is a regular groupoid of type zero. Moreover, since the core is trivial, $T\G\cong s^*TM$ is given by the transformation groupoid of a representation of $\G$ on $TM$. It follows from Example \ref{example: Sec of repns} that the space of multiplicative vector fields $\mathfrak{X}_{mult}(\G)$ is isomorphic to the space $\Gamma(TM)^{\G}$ of $\G$-invariant vector fields on $M$. 

Let us see a particular example of \'etale groupoid.  Assume that $\mathcal{F}$ is a regular foliation on $M$ and let $\mathrm{Hol}(\mathcal{F})\rr M$ be the holonomy groupoid. If $S\subseteq M$ is a complete transversal of $\mathcal{F}$, i.e. $S$ intersects transversally each leaf at least once, then the inclusion $i_S:S\rmap M$ can be used to construct the \textbf{restricted holonomy groupoid} $\mathrm{Hol}_S(\mathcal{F}):=i^!_S\mathrm{Hol}(\mathcal{F})\rr S$ which is just the pullback groupoid by the inclusion map. The canonical map $\mathrm{Hol}_S(\mathcal{F})\rmap \mathrm{Hol}(\mathcal{F})$ is a Morita map, hence the complex of multiplicative vector fields on $\mathrm{Hol}(\mathcal{F})\rr M$ is quasi-isomorphic to the complex of multiplicative vector fields on  $\mathrm{Hol}_S(\mathcal{F})\rr S$. The latter is an \'etale groupoid, so we conclude that its complex of multiplicative vector fields is

$$0\rmap \Gamma(TS)^{\mathcal{F}},$$

\noindent where $\Gamma(TS)^{\mathcal{F}}$ is the space of vector fields on $S$ which are invariant by the holonomy representantion.

\end{example}

\begin{example}(Regular groupoids)\label{vf:regular}

If $\G\rr M$ is a regular Lie groupoid, then the tangent groupoid $T\G\rr TM$ is a regular $\vb$-groupoid, whose complex of multiplicative sections has been described in Example \ref{ex:complexregularTG}. Due to Proposition \ref{regularqisotypezero}, $T\G$ is $\vb$-Morita equivalent to its type zero component $t^*K\oplus s^*\nu\rr \nu$, where $K=\ker(\rho)$ and $\nu=TM/\mathrm{im}(\rho)$. In particular, the complex of multiplcative vector fields is quasi-isomorphic to the 2-term complex $\Gamma(K)\rmap \Gamma(s^*\nu)$.

A particular class of examples of regular groupoids is given by the so-called \emph{foliation groupoids}, i.e. Lie groupoids which are Morita equivalent to an \'etale groupoid. It was shown in \cite{crainicm} that a Lie groupoid $\G$ with Lie algebroid $A$ is a foliation groupoid if and only if the anchor map $\rho:A\rmap TM$ is injective. As as result, the complex of multiplicative vector fields on a foliation groupoid is quasi-isomorphic to $0\rmap \Gamma(s^*\nu)$.

\end{example}

\subsection{Vector fields on differentiable stacks}
\label{subsection: vect on stacks}

We start by a quick review on differentiable stacks and their relation with Lie groupoids. Let $Man$ be the site of smooth manifolds and $Grpd$ the category of groupoids. A \textbf{stack} over manifolds is a fibered category $\pi:\mathfrak{S}\rmap Man$ satisfying a descent condition \cite{bx, vistoli}. Given a smooth manifold $M$, the fiber $\mathfrak{S}(M):=\pi^{-1}(M)$ has a natural structure of groupoid. Moreover, the choice of a \emph{cleavage} of $\pi:\mathfrak{S}\rmap Man$ defines a pseudo-functor $\mathfrak{S}:Man^{op}\rmap Grpd$ which satisfies a gluing condition at both the level of objects and morphisms. This establishes a correspondence between fibered categories together with a \emph{cleavage} and pseudo-functors $Man^{op}\rmap Grpd$ defining a sheaf of groupoids. For more details see \cite{bx, vistoli}.

\begin{example}(Manifolds as stacks)\label{ex:stackmanifold}

Every smooth manifold $M$ can be seen as a stack $[M]\rmap Man$, where $[M]$ is the category whose objects are smooth maps $U\rmap M$ and morphisms are commutative triangles. 

\end{example}

\begin{example}(Lie groups as stacks)\label{ex:classifyingstack}

Let $G$ be a Lie group. The \textbf{classifying stack} of $G$ is defined as the category $\pi:[*/G]\rmap Man$, where $[*/G]$ is the category with objects being principal $G$-bundles and morphisms are principal $G$-bundle maps. This defines a stack $[M/\G]\rmap Man; P\to N\mapsto N$, referred to as the \textbf{quotient stack} of $\G\rr M$.

\end{example}

\begin{example}(Quotient stack of a Lie group action)\label{ex:stackaction}

Consider a Lie group $G$ acting on a smooth manifold $M$. The \textbf{quotient stack} of the action is the fibered category $[M/G]$ whose objects are pairs $(P\to N, P\to M)$ where $P\rmap N$ is a principal $G$-bundle and $P\rmap M$ is a $G$-equivariant smooth map. Morphisms $(P\to N, P\to M)\rmap (P'\to N', P'\to M)$ are given by principal $G$-bundle maps $\varphi:P\rmap P'$ which determine a commutative triangle with target $M$.

\end{example}

All the examples above can be unified by looking at the quotient stack of a Lie groupoid, which we describe below.

\begin{example}(Quotient stack of a Lie groupoid)\label{ex:quotientstack}

Let $\G\rr M$ be a Lie groupoid. Consider the fibered category $[M/\G]$ of principal $\G$-bundles and morphisms of principal $\G$-bundles.

\end{example}

A \textbf{morphism} of stacks is just a morphism between the underlying fibered categories. An \textbf{atlas} of a stack $\mathfrak{S}$ is a \emph{representable epimorphism} of stacks $[M]\rmap \mathfrak{S}$, where $[M]$ is the stack associated to a manifold. That is, a stack morphism satisfying the property that for any other map $[N]\rmap \mathfrak{S}$, the fibered product of stacks $[M]\times_{\mathfrak{S}}[N]$ is a manifold and the projection $[M]\times_{\mathfrak{S}}[N]\rmap [N]$ is a surjective submersion. A stack $\mathfrak{S}$ is called \textbf{differentiable} if it admits an atlas. 

It is well known that if $[M]\rmap \mathfrak{S}$ is an atlas of a stack $\mathfrak{S}$, then the fibered product $[M]\times_{\mathfrak{S}}[M]$ has a canonical structure of Lie groupoid over $M$, with source and target maps given the canonical projections. Moreover, the stack $\mathfrak{S}$ is isomorphic to the quotient stack of the Lie groupoid $[M]\times_{\mathfrak{S}}[M]\rr M$. Conversely, given a Lie groupoid $\G\rr M$, the canonical map $[M]\rmap [M/\G]; (\phi:N\rmap M)\mapsto t\circ \phi:\phi^*\G\rmap N$, defines an atlas of $[M/\G]$ and the associated Lie groupoid $[M]\times_{[M/\G]}[M]\rr M$ is Morita equivalent to $\G\rr M$.

Given a differentiable stack $\mathcal{S}$, a Lie groupoid $\G\rr M$ with $[M/\G]\cong \mathfrak{S}$ is referred to as a \textbf{presentation} of $\mathfrak{S}$.

\begin{remark}
In this section, we study differentiable stacks via Lie groupoids. We think of a differentiable stack in two different ways: \begin{enumerate}
\item as the quotient stack $[M/\G]$ of a Lie groupoid $\G\rr M$, and
\item as a \emph{generalized quotient} of smooth manifolds, i.e. a representable epimorphism or atlas $[M]\rmap \mathfrak{S}$.
\end{enumerate}

\end{remark}

\begin{remark}
In \cite{hepworth} Hepworth defines the tangent stack $T\mathcal X$ of a stack $\mathcal X$, and the category of sections $\Gamma(\mathcal X,T\mathcal X)$. If $\mathcal X \simeq [M/\mathcal G]$ then it is shown in \cite{hepworth} that there is an equivalence of categories
\begin{equation}
\label{eqn: hepworth equivalence}
\mathrm{Sec}(\mathcal G,T\mathcal G)
\simeq 
\Gamma( \mathcal X, T\mathcal X) 
\end{equation}
which identifies the zero sections in each category. Recall from Section \ref{subsection: 2-vector spaces} that if $V^\bullet$ is a 2-term complex then $H^0 (V^\bullet)$ can be identified with the group $\pi_1 \left( V_1 \rightrightarrows V_0\right)$ of automorphisms of the zero object in the corresponding 2-vector space $V_1 \rightrightarrows V_0$, and $H^1(V^\bullet)$ can be identified with the set $\pi_0 \left( V_1 \rightrightarrows V_0\right)$ of isomorphism classes of objects in  $V_1 \rightrightarrows V_0$. Therefore, via the equivalence of categories of Theorem \ref{2vect} the equivalence (\ref{eqn: hepworth equivalence}) determines bijections
\begin{align}
\label{eqn: H0 bijection}
H^0 (L_\mathcal G) &  \cong
\pi_1 \Gamma(\mathcal X,T \mathcal X) \\
\label{eqn: H1 bijection}
H^1 (L_\mathcal G) &  \cong
\pi_0 \Gamma(\mathcal X,T \mathcal X)
\end{align}
 In particular, by Proposition \ref{vf:moritainvariance} the sets $\pi_1 \Gamma(\mathcal X,T \mathcal X)$ and $\pi_0 \Gamma(\mathcal X,T \mathcal X)$ carry Lie algebra structures, well defined up to isomorphism. (Note that the mismatch between the gradings on the left and right hand sides of (\ref{eqn: H0 bijection}) and (\ref{eqn: H1 bijection}) arises from our conventions on 2-vector spaces.)
\end{remark}

Our main goal now is to introduce a geometric notion of vector field on a differentiable stack. For that, we first study vector fields on \emph{quotient manifolds}, that is, a smooth manifold $N$ which is the base of a surjective submersion $\pi:M\rmap N$. 

\begin{remark}(Vector fields on quotient manifolds)\label{rmk:vfquotientmanifold}

It is well known that there is a canonical identification between vector fields on $N$ and horizontal vector fields on $M$, that is

\begin{equation}\label{eq:identificationvfquotientmanifold}
\mathfrak{X}(N)\cong \frac{\mathfrak{X}(M)^{\pi}}{\Gamma(\ker(T\pi))},
\end{equation}

\noindent where $\mathfrak{X}(M)^{\pi}$ denotes the space of $\pi$-projectable vector fields. The surjective submersion $\pi:M\rmap N$ can be seen as an atlas of $N$ viewed as a differentiable stack. Hence, the submersion groupoid $M\times_NM\rr M$ is a presentation of the stack $[N]$. As observed in Example \ref{vf:submersions}, the complex of multiplicative vector fields on the submersion groupoid $M\times_NM\rr M$ is given by $\Gamma(\ker(T\pi))\rmap \mathfrak{X}(M)^{\pi}$ and the first cohomology group is $\frac{\mathfrak{X}(M)^{\pi}}{\Gamma(\ker(T\pi))}$, which identifies with the space of vector fields on $N$ via \eqref{eq:identificationvfquotientmanifold}. 
\end{remark}

This motivates the following definition.

\begin{definition}\label{vfstacks}

Let $\G\rr M$ be a Lie groupoid. A \textbf{geometric vector field} on the quotient stack $[M/\G]$ is a degree one cohomology class $[V]\in H^1(L_{\G})$ of a multiplicative vector field $V\in\mathfrak{X}_{mult}(\G)$.
\end{definition}

The space of vector fields on the quotient stack $[M/\G]$ is 

$$\mathfrak{X}([M/ \G]):=H^1(L_{\G}).$$
According to Corollary \ref{cohomologyliealgebra}, the space of vector fields $\mathfrak{X}([M/\G])$ inherits a natural Lie algebra structure. Moreover, as a consequence of Theorem \ref{theorem: invariance of Lie 2-algebras}, if $\G\rr M$ and $\G'\rr M'$ are Morita equivalent Lie groupoids, then the corresponding Lie algebras of vector fields $\mathfrak{X}([M/\G])$ and $\mathfrak{X}([M'/\G'])$ are isomorphic.

\begin{remark}(Stacks as generalized quotients of manifolds)

This corresponds to the stacky version of the description of vector fields on quotient manifolds explained in Remark \ref{rmk:vfquotientmanifold}. Let $\mathfrak{S}$ be a differentiable stack. If $\pi:[M]\rmap \mathfrak{S}$ is an atlas, then the Lie groupoid $[M]\times_{\mathfrak{S}}[M]\rr M$ can be thought of as kind of \emph{submersion groupoid} associated to the representable epimorphism $\pi:[M]\rmap \mathfrak{S}$. In this case, the Lie groupoid $[M]\times_{\mathfrak{S}}[M]\rr M$ is a presentation of $\mathfrak{S}$ and hence, the space of geometric vector fields on $\mathfrak{S}$ is isomorphic to $H^1(L_{[M]\times_{\mathfrak{S}}[M]})$. Notice that if $[M']\rmap \mathfrak{S}$ is another atlas, then $[M']\times_{\mathfrak{S}}[M']\rr M'$ is Morita equivalent to $[M]\times_{\mathfrak{S}}[M]\rr M$ and hence $H^1(L_{[M']\times_{\mathfrak{S}}[M']})\cong H^1(L_{[M]\times_{\mathfrak{S}}[M]})$ due to Theorem \ref{moritainvarianceVB}.
\end{remark}

In general, the space of geometric vector fields on a quotient stack $[M/\G]$ can be described by using representations up to homotopy. More precisely, if one applies Corollary \ref{cohomology} to the special case of multiplcative vector fields, the degree one cohomology of the complex of multiplicative vector fields on a Lie groupoid $\G\rr M$ coincides with $H^1(\G,Ad_\G)$, the cohomology of $\G$ with coefficients in the adjoint representation up to homotopy. However, describing $H^1(\G,Ad_\G)$ is not easy in general. In what follows, we describe the Lie algebra of geometric vector fields of some particular differentiable stacks.

\begin{example}(Manifolds)

The stack $[M]$ associated to a smooth manifold is presented by the unit groupoid $M\rr M$. The description of the complex of multiplicative vector fields on $M\rr M$ given in Example \ref{vf:manifolds} implies that $\mathfrak{X}([M])\cong \mathfrak{X}(M)$.

\end{example}

\begin{example}(The classifying stack of a Lie group)
\label{example: classifying stack of a Lie group}
Let $[*/G]$ be the classifying stack associated to a Lie group described in Example \ref{ex:classifyingstack}. The Lie groupoid $G\rr \{*\}$ is a presentation of $[*/G]$ and by Example \ref{vf:liegroup}, the space of vector fields on the classifying stack is given by $\mathfrak{X}([*/G])\cong H^1(G,\mathfrak{g})$ the first cohomology group of $G$ with coefficients in the adjoint representation.

\end{example}

\begin{example}(The quotient stack of a Lie group action)
\label{example: quotient stack of a Lie group action}

Let $[M/G]$ be the quotient stack of the action of a Lie group on a smooth manifold $M$. A presentation of $[M/G]$ is given by the transformation groupoid $G\ltimes M\rr M$. Hence, the space of geometric vector fields on $[M/G]$ is given by $H^1(L_{G\ltimes M})$. Let us discuss some special cases of this situation:

\begin{itemize}

\item[i)] If $G$ is compact, a result of Hepworth \cite{hepworth} guarantees that the complex of multiplicative vector fields on $G\ltimes M$ is quasi-isomorphic to the complex $C^{\infty}(M,\mathfrak{g})^G\rmap \mathfrak{X}(M)^G$ (see Example \ref{vf:liegroupaction}). As a consequence, the space of geometric vector fields on $[M/G]$ identifies with certain quotient of the space of invariant vector fields $\mathfrak{X}(M)^G$. For instance, if $G$ is a finite group, then the space of geometric vector fields on the quotient stack $[M/G]$ identifies with the space of invariant vector fields $\mathfrak{X}(M)^G$.

\item[ii)] If $G$ acts on $M$ in a proper and locally free manner, then the transformation groupoid $G\ltimes M\rr M$ is regular. Let $\nu=TM/\mathrm{im}(\rho)$ the transversal bundle of the action as in Example \ref{vf:locallyfree}. In this case, the transformation groupoid is regular and the space of geometric vector fields on the quotient stack $[M/G]$ identifies with the space $\Gamma(\nu)^G$ of transversal $G$-invariant vector fields. 

\end{itemize}

\end{example}

\begin{example}(Orbifolds)

It is well known that orbifolds can be seen as certain class of differentiable stacks. More precisely, orbifolds are presented by proper and \'etale Lie groupoids. Assume that $\G\rr M$ is a proper and \'etale Lie groupoid. The Lie algebroid of $\G$ is $A=M\times 0$ and the complex of multiplicative vector fields on $\G$ is $0\rmap \mathfrak{X}(M)^\G$ (see Example \ref{vf:etale}). Therefore, the space of geometric vector fields on the quotient stack/orbifold $[M/\G]$ is given by $\mathfrak{X}([M/\G])=\mathfrak{X}(M)^\G$. A particular instance of this situation is that of the restricted holonomy groupoid. Consider $\mathrm{Hol}(\mathcal{F})\rr M$ the holonomy groupoid of a regular foliation on $M$. If $S\subset M$ is a complete transversal, then the restricted holonomy groupoid $\mathrm{Hol}_S(\mathcal{F})\rr S$ is \'etale

\end{example}

\begin{example}(Quotient stack of regular groupoids)
\label{example: Quotient stack of regular groupoids)}
Given a regular Lie groupoid $\G\rr M$, the tangent groupoid $T\G\rr TM$ is a regular $\vb$-groupoid. Therefore, one can describe geometric vector fields on the corresponding quotient stack $[M/\G]$ in terms of the type zero part of $T\G$, i.e. via transversal information. It follows from the discussion in Example \ref{vf:regular} that the space of geometric vector fields on the stack $[M/\G]$ identifies with the first cohomology group of the complex $\Gamma(K)\rmap \Gamma(s^*\nu)$ where $K=\ker(\rho)$ and $\nu=TM/\mathrm{im}(\rho)$.

A particular instance of this situation is given by foliation groupoids. If $\G\rr M$ is a foliation groupoid (its Lie algebroid has injective anchor map), then the space of geometric vector fields on $[M/\G]$ identifies with $\Gamma(s^*\nu)$. Also, vector fields on the stack associated to the holonomy groupoid of a foliation can be described as follows. Consider $\mathrm{Hol}(\mathcal{F})\rr M$ the holonomy groupoid of a regular foliation on $M$. If $S\subset M$ is a complete transversal, then the restricted holonomy groupoid $\mathrm{Hol}_S(\mathcal{F})\rr S$ is \'etale and hence $\mathfrak{X}([M/\mathcal{F}])\cong \mathfrak{X}(S)^{\mathcal{F}}$.

\end{example}

It is interesting to study in which sense the category of multiplicative vector fields and geometric vector fields on quotient stacks can be seen as derivations of a suitable notion of algebra of functions on differentiable stacks. This will be treated in a future work \cite{OW}.

\subsection{The deformation complex}

In \cite{CMS} the authors define and study the \emph{deformation complex} $C^\bullet_{def}(\mathcal G)$ of a Lie groupoid $\mathcal G$. Propositions 3.9 and 4.3 in loc. cit. show that 
\[
C^\bullet_{mult}\left(T\mathcal G\right) \cong
C^\bullet_{def}\left(\mathcal G\right)^{\le 1}
\]
where $C^\bullet_{def}\left(\mathcal G\right)^{\le 1}$ is the 2-term truncation of $C^\bullet_{def}\left(\mathcal G\right)$. Some of the results in \cite{CMS} then provide alternative proofs of some of computations in the examples in Section \ref{subsection: vect on stacks} above. In particular, Proposition 3.1 in loc. cit. gives an alternative proof of the computation in Example \ref{example: classifying stack of a Lie group} above, Proposition 3.3 of loc. cit. of that of Example \ref{example: Quotient stack of regular groupoids)}, and Theorem 6.1 in loc. cit. of that of the discussion of free proper actions in Example \ref{example: quotient stack of a Lie group action}.

\appendix

\section{Appendix}

\subsection{Fibre-wise linear functions}

If $E \to M$ is a vector bundle then we denote by $l$ the canonical isomorphism
\begin{align*}
l : \Gamma(E) & \to C^\infty_{lin}\left(E^\ast\right) \\
e & \mapsto \left( \alpha \in E^\ast_m \mapsto \alpha \left( e(m)\right) \right) 
\end{align*}
where $C^\infty_{lin}(E^\ast)$ is the vector space of fibre-wise linear smooth functions on $E^\ast$.

\begin{proposition}
\label{prop: naturality of l}
\begin{enumerate}
\item
If $\Phi : E \to F$ is a morphism of vector bundles over $M$
 then the following diagram commutes:
\[
\xymatrix{
C^\infty_{lin}\left(E^\ast\right)  \ar[r]^{\left({\Phi^\ast}\right)^\ast} & 
C^\infty_{lin}\left(F^\ast\right) \\
\Gamma(E)\ar[u]^l \ar[r]_{\Phi} & \Gamma(F) \ar[u]_l
}
\]
(Here, $\Phi^\ast : F^\ast \to E^\ast$ is the morphism of vector bundles over $M$ dual to $\Phi$, and $\left(\Phi^\ast\right)^\ast$ is the pullback of smooth functions $\left(\Phi^\ast\right)^\ast : C^\infty\left(E^\ast\right) \to C^\infty\left(F^\ast\right)$ induced by $\Phi^\ast$. Note that $\left(\Phi^\ast\right)^\ast$ restricts to a map $\left(\Phi^\ast\right)^\ast : C^\infty_{lin}\left(E^\ast\right) \to C^\infty_{lin}\left(F^\ast\right)$ because $\Phi^\ast$ is fibre-wise linear.)
\item
If $E \to M$ is a vector bundle, $f: N \to M$ is a smooth map, and $f^! : f^\ast E \to E$ is the canonical pullback morphism covering $f$
\[
\xymatrix{
f^\ast E  \ar[d]  \ar[r]^{f^!}  &   E \ar[d] \\
N  \ar[r]_f & M
}
\]
then the following diagram commutes:
\[
\xymatrix{
C^\infty_{lin} \left(f^\ast E^\ast\right)  &  \ar[l]_-{\left({f^!}\right)^\ast}  C^\infty_{lin} \left(E^\ast \right)  \\
\Gamma\left( f^\ast E \right) \ar[u]^l  & \ar[l]^{f^\ast}  \Gamma(E) \ar[u]_l
}
\]
\end{enumerate}
\end{proposition}

\begin{proof}
If $e \in \Gamma(E)$, $e' \in \Gamma(F)$, $x \in N$, $\alpha \in E^\ast_{f(x)}$, and $\bar \alpha = \left(x,\alpha\right) \in \left(f^\ast E^\ast\right)_x$ , then
\begin{align*}
\left(\left( \Phi^\ast\right)^\ast l_e \right) (\bar\alpha)  & =  l_e \left( \Phi^\ast \bar\alpha \right) \\
& = \left( \Phi^\ast \bar\alpha\right) \left(e(x)\right) \\
& = \bar\alpha \left(\Phi \left( e(x)\right) \right) \\
& = l_{ \Phi(e)} (\bar\alpha) \\
\\
\left(\left(f^!\right)^\ast  l_{e}\right) \left(x,\alpha\right) 
& = l_e  \left(f^! (x,\alpha) \right) \\
& = l_e \left(\alpha\right) \\
& = \alpha \left( e(f(x)\right) \\
& = \left(x,\alpha\right) \left(x,e(\phi(x))\right) \\
& = l_{f^\ast e}\left(g,\alpha\right)
\end{align*}
This shows that 
\begin{align*}
\left(\Phi^\ast \right)^\ast \circ l & = l \circ  \Phi \\
\left(f^!\right)^\ast \circ l & = l \circ f^\ast
\end{align*}
as required.
\end{proof}

\subsection{Weak fibre products of $\mathcal{VB}$ and $\mathcal{LA}$-groupoids}

The purpose of this appendix is to show the existence of weak fibre products of $\mathcal{VB}$ and $\mathcal{LA}$ groupoids under certain conditions. This is needed is Sections \ref{sec:moritainvariance}  and \ref{section: Morita equivalence LA}. We start by recalling the construction of weak fibre products of Lie groupoids. See \cite{MoerdijkMrcun} for more details.

Let $\mathcal G \rightrightarrows M,\mathcal G' \rightrightarrows M'$ and $\mathcal H\rightrightarrows N$ be Lie groupoids, with structure maps denoted $s,s'$ and $\tilde s$ and similar. Let $\phi:\mathcal G \to \mathcal H$ and $\phi':\mathcal G' \to \mathcal H$ be morphisms of groupoids as in the following diagram:
\begin{equation}
\label{diag: pullback diagram for weak pullback of groupoids}
\xymatrix{
&  \mathcal G' \ar[d]^{\phi'} \\
\mathcal G \ar[r]_{\phi}  &  \mathcal H
}
\end{equation}
Consider the following groupoid $\mathcal G \times_\mathcal H^w \mathcal G'$ (in sets):
\begin{itemize}
\item  The objects are triples $(x,h,x')$, where $x \in \mathcal M, x' \in M'$, and $h \in \mathcal H\left(\phi_0(x),\phi_0'(x')\right)$.
\item A morphism $(x,h,x') \to (y,k,y')$ is given by a pair $(g,g') \in \mathcal G\left(x,y\right) \times \mathcal G'\left(x',y'\right)$, such that $\phi'(g')h = k\phi(g)$.
\item The structure maps are determined component-wise by those of $\mathcal G$ and $\mathcal G'$.
\end{itemize}
The sets of objects and morphisms can be identified with certain fibre products over $N$:
\begin{align}
\label{eqn: weak fibre product morphisms}
\mathcal G \times_\mathcal H^w \mathcal G' & = 
\mathcal G \times_N^{\phi_0s,\tilde s} \mathcal H \times_N^{\tilde t,\phi_0's'} \mathcal G' \\
\left(\mathcal G \times_\mathcal H^w \mathcal G'\right)_0 & = 
M \times_N^{\phi_0,\tilde s} \mathcal H \times_N^{\tilde t,\phi_0'} M'
\label{eqn: weak fibre product objects}
\end{align}
If the fibre products of smooth manifolds (\ref{eqn: weak fibre product morphisms}) and (\ref{eqn: weak fibre product objects}) exist then $\mathcal G \times_\mathcal H^w \mathcal G'$ is a Lie groupoid  and the canonical projections $\mathcal G\times_\mathcal H^w \mathcal G' \to \mathcal G$ and $\mathcal G\times_\mathcal H^w \mathcal G'$ are morphisms of Lie groupoids. If $\phi$ is a weak equivalence then so is the projection $\mathcal G\times_\mathcal H^w \mathcal G' \to \mathcal G'$.

\begin{proposition}
\label{prop: existence of VB weak fibre products} 
Let $\mathcal V \rightrightarrows E$, $\mathcal V' \rightrightarrows E'$, $\mathcal W \rightrightarrows F$ be $\mathcal{VB}$-groupoids with bases $\mathcal G \rightrightarrows M$, $\mathcal G' \rightrightarrows M'$, $\mathcal H \rightrightarrows N$. Let $\Phi:\mathcal V \to \mathcal W$ and $\Phi':\mathcal V' \to \mathcal W$ be morphisms of $\mathcal{VB}$ groupoids covering morphisms of Lie groupoids $\phi:\mathcal G \to \mathcal H$ and $\phi':\mathcal G' \to \mathcal H$ as in the following diagrams:
\[
\xymatrix{
&  \mathcal V'  \ar[d]^{\Phi'}  &  &
&    \mathcal G' \ar[d]^{\phi'}
\\
\mathcal V  \ar[r]_{\Phi}  &   \mathcal W &  &
  \mathcal G  \ar[r]_{\phi}  &  \mathcal H
}
\]
Suppose that $\Phi$ is a $\mathcal{VB}$-Morita morphism (Definition \ref{def: VB Morita maps}). Then:
\begin{enumerate}
\item The weak fibre products of Lie groupoids $\mathcal V\times_\mathcal W^w \mathcal V'$ and $\mathcal G \times_\mathcal H ^w \mathcal G'$ exist.
\item The Lie groupoid $\mathcal V\times_\mathcal W^w \mathcal V'$ is a $\mathcal{VB}$-groupoid over $\mathcal G \times_\mathcal H ^w \mathcal G'$ in a canonical way:
\[
\xymatrix{
\mathcal V \times_\mathcal W^w \mathcal V'  \ar[d]
\ar@<.5ex>[r]\ar@<-.5ex>[r] & \left(\mathcal V \times_\mathcal W \mathcal V'\right)_0  \ar[d]\\
\mathcal G \times_\mathcal H^w \mathcal G'   
\ar@<.5ex>[r]\ar@<-.5ex>[r] & \left(\mathcal G \times_\mathcal H \mathcal G'\right)_0
}
\]
\item  The projections $\mathcal V\times_\mathcal W^w \mathcal V' \to \mathcal V$ and $\mathcal V\times_\mathcal W^w \mathcal V' \to \mathcal V'$ are morphisms of $\mathcal{VB}$-groupoids.
\item The projection $\mathcal V\times_\mathcal W^w \mathcal V' \to \mathcal V'$ is a $\mathcal{VB}$-Morita morphism.
\item As morphisms of Lie groupoids the projections $\mathcal V \times_\mathcal W^w \mathcal V' \to \mathcal V'$ and $\mathcal G \times_\mathcal H^w \mathcal G' \to \mathcal G'$ are surjective on objects, morphisms and composable pairs.
\end{enumerate}
\end{proposition}

\begin{proof}
Consider the following diagrams of morphisms of vector bundles:
\begin{equation}
\label{diag: weak fibre product 1}
\xymatrix{
\mathcal V \ar[dd] \ar[dr]^{\Phi_0 \tilde s} & & \mathcal W \ar[dd] \ar[dl]_{\tilde s} \ar[dr]^{\tilde t} & & \mathcal V' \ar[dd] \ar[dl]_{\Phi'_0 \tilde s'} \\
& F \ar[dd] & & F \ar[dd] \\
\mathcal G \ar[dr]_{\phi_0 s} & &  \ar[dl]^{s} \mathcal H \ar[dr]_{t} &  & \mathcal G' \ar[dl]^{\phi'_0 s'} \\
& N &  &  N
}   
\end{equation}

\begin{equation}
\label{diag: weak fibre product 0}
\xymatrix{
E \ar[dd] \ar[dr]^{\Phi_0} & & \mathcal W \ar[dd] \ar[dl]_{\tilde s} \ar[dr]^{\tilde t} & & E' \ar[dd] \ar[dl]_{\Phi'_0} \\
& F \ar[dd] & & F \ar[dd] \\
M \ar[dr]_{\phi_0} & &  \ar[dl]^{s} \mathcal H \ar[dr]_{t} &  & M' \ar[dl]^{\phi'_0} \\
& N &  &  N 
}
\end{equation}
We first show that the iterated fibre products of the diagrams (\ref{diag: weak fibre product 1}) and (\ref{diag: weak fibre product 0}) exist in the category of vector bundles. Note that the `base component' of a morphism of vector bundles is a submersion whenever the `total space component' is, and that transverse fibre products of vector bundles always exist. Therefore, in order to show that a certain fibre product of vector bundles exists it is sufficient to show that the total space component of one of the morphisms is a submersion. With this in mind, the proof is identical to that of the analogous result for Lie groupoids given in \cite{MoerdijkMrcun}, but we spell out the details for convenience.

We start with the diagram (\ref{diag: weak fibre product 0}). As $\tilde s$ is a submersion the fibre product
\[
\xymatrix{
E \times_F ^{\Phi_0,\tilde s} \mathcal W \ar[d] \\
M \times_N^{\phi_0,s} \mathcal H
}
\]
exists. By assumption $\Phi$ is a $\mathcal{VB}$-Morita morphism, and so $\Phi$ is in particular a weak equivalence of Lie groupoids. Therefore the map
\begin{align*}
\tilde t \mathrm{pr}_2 : E \times_F ^{\Phi_0,\tilde s} \mathcal W & \to  \mathcal W 
\end{align*}
is a submersion. It follows that the fibre product
\begin{equation}
\label{diag: fibre product bundle 0}
\xymatrix{
E \times_F^{\Phi_0,\tilde s} \mathcal W \times_F^{\tilde t,\Phi_0'} E'
\ar[d] \\
M \times_N^{\phi_0,s} \mathcal H \times_N^{t,\phi_0'} M'
}
\end{equation}
of the diagram (\ref{diag: weak fibre product 0}) exists. This implies that the fibre product
\begin{equation}
\label{diag: fibre product bundle 1}
\xymatrix{
\mathcal V \times_F^{\Phi_0 \tilde s,\tilde s} \mathcal W \times_F^{\tilde t,\Phi_0' \tilde s'} \mathcal V'
\ar[d] \\
\mathcal G \times_N^{\phi_0s,s} \mathcal H \times_N^{t,\phi_0's'} \mathcal G'
}
\end{equation}
of the diagram (\ref{diag: weak fibre product 1}) exists because it can be constructed as the iterated fibre product
\[
\xymatrix{
\mathcal V \times_F^{\Phi \tilde s,\tilde s}\mathcal W\times_F^{\tilde t, \tilde s' \Phi_0'} \mathcal V'   \ar[d]  \ar[rr]  & & 
\mathcal V' \ar[d]^{\tilde s'}
\\
\mathcal V \times_F^{\Phi_0 \tilde s}\mathcal W\times_F^{\tilde t, \Phi_0} E' \ar[d]_{\mathrm{pr}_1} \ar[r] &  
E \times_F^{\Phi_0,\tilde s} \mathcal W\times_F^{\tilde t,\Phi_0'} E'  \ar[d]^{\mathrm{pr}_1} \ar[r]_-{\mathrm{pr}_3} & E'  \\
\mathcal V \ar[r]_{\tilde s}  &  E
}
\]
which exists because $\tilde s$ and $\tilde s'$ are submersions. Therefore, the weak fibre products $\mathcal V\times_\mathcal W^w \mathcal V'$ and $\mathcal G \times_\mathcal H ^w \mathcal G'$ exist, proving Statement 1, and there are canonical vector bundle structures $\mathcal V\times_\mathcal W^w \mathcal V' \to \mathcal G \times_\mathcal H ^w \mathcal G'$ and  $(\mathcal V\times_\mathcal W^w \mathcal V')_0 \to (\mathcal G \times_\mathcal H ^w \mathcal G')_0$.

Statements 2 and 3 follow from the fact that the structure maps of the Lie groupoid $\mathcal V \times_\mathcal W^w\mathcal V'$ can be expressed in terms of fibre products and compositions of the structure maps of $\mathcal V,\mathcal W$ and $\mathcal V'$, and the maps $\Phi$ and $\Phi'$, all of which are vector bundle morphisms covering the corresponding maps between $\mathcal G,\mathcal H$ and $\mathcal G'$. Using the fact that $\Phi:\mathcal V \to \mathcal W$ is a $\mathcal{VB}$-Morita morphism and therefore a weak equivalence of Lie groupoids,  Statement 4, and the fact that the projections $\mathcal V \times_\mathcal W^w \mathcal V' \to \mathcal V'$ and $\mathcal G\times_\mathcal H^w\mathcal G' \to \mathcal G'$ are surjective on objects follows from Proposition 5.12 of \cite{MoerdijkMrcun}.
Statement 5 follows from the fact that a weak equivalence which is surjective on objects is also surjective on morphisms and composable pairs.
\end{proof}

\begin{proposition}
\label{prop: existence of LA weak fibre products} 
Proposition \emph{\ref{prop: existence of VB weak fibre products}} continues to hold if $\mathcal{VB}$-groupoids are replaced by $\mathcal{LA}$-groupoids throughout.
\end{proposition}

\begin{proof}
The proof is identical to that of Proposition \ref{prop: existence of VB weak fibre products}: the fibre products which appear in the proof are all transverse fibre products of Lie algebroids - which always exist (see, 
\cite{HigginsMackenzie},\cite{Stefanini},\cite{BursztynCabreradelHoyo}).
\end{proof}

\end{document}